\documentclass[a4paper, 11pt, twoside, openany]{article} 
\usepackage{amsmath,amssymb,stmaryrd, mathtools}
\usepackage{fancyhdr}
\usepackage{footmisc}
\usepackage{hyperref}
\usepackage[margin=22mm]{geometry}
\usepackage[UKenglish]{datetime}
\usepackage{graphicx}
\usepackage[usenames, dvipsnames]{color}
\usepackage{accents}
\usepackage[all]{xy}

\numberwithin{equation}{section}

\newcounter{savefootnote}
\newcounter{symfootnote}
\newcommand{\symfootnote}[1]{%
   \setcounter{savefootnote}{\value{footnote}}%
   \setcounter{footnote}{\value{symfootnote}}%
   \ifnum\value{footnote}>8\setcounter{footnote}{0}\fi%
   \let\oldthefootnote=\thefootnote%
   \renewcommand{\thefootnote}{\fnsymbol{footnote}}%
   \footnote{#1}%
   \let\thefootnote=\oldthefootnote%
   \setcounter{symfootnote}{\value{footnote}}%
   \setcounter{footnote}{\value{savefootnote}}%
}

\newcommand{\pp}[2]{\frac{\partial #1}{\partial #2}} 
\newcommand{\dd}[2]{\frac{\delta #1}{\delta #2}}

\newtheorem{definition}{Definition}
\newtheorem{prop}{Proposition}
\newtheorem{lemma}{Lemma}
\newtheorem{remark}{Remark}
\begin{document}
\begin{center}
{\Large Energy conserving SUPG methods for compatible finite element schemes in numerical weather prediction}
\end{center}
\vspace{-3mm}
\hrulefill
\begin{center}
{Golo A. Wimmer$^{1}$\symfootnote{Current affiliation: Los Alamos National Laboratory; correspondence to: gwimmer@lanl.gov}, Colin J. Cotter$^1$ and Werner Bauer$^{1,2}$}\\
\vspace{2mm}
{\textit{$^1$Imperial College London, $^2$INRIA Rennes}}\\
\vspace{4mm}
\today
\end{center}
\vspace{5mm}
\begin{center}
\textbf{Abstract}
\end{center}
We present an energy conserving space discretisation based on a Poisson bracket that can be used to derive the dry compressible Euler as well as thermal shallow water equations. It is formulated using the compatible finite element method, and extends the incorporation of upwinding for the shallow water equations as described in Wimmer, Cotter, and Bauer (2020). While the former is restricted to DG upwinding, an energy conserving SUPG method for the (partially) continuous Galerkin thermal field space is newly introduced here. The energy conserving property is validated by coupling the Poisson bracket based spatial discretisation to an energy conserving time discretisation. Further, the discretisation is demonstrated to lead to an improved thermal field development with respect to stability when upwinding is included. An approximately energy conserving scheme that includes upwinding for all prognostic fields with a smaller computational cost is also presented. In a falling bubble test case used for the Euler equations, the latter scheme is shown to resolve small scale features at coarser resolutions than a corresponding scheme derived directly from the equations without the Poisson bracket framework.\\ \\
%
%
%
%
%
\section{Introduction}
Finite element methods have recently gained an increased interest in numerical weather prediction (NWP), as they allow for higher order discretisations and more general meshes, thus avoiding the parallel computing issues associated with grid poles. This includes spectral element methods, discontinuous Galerkin methods, and the compatible finite element method \cite{COTTER20127076}, where in the latter finite element spaces are mapped to one another via differential operators \cite{natale2016compatible}. In the context of NWP, the compatible finite element method can be seen as an extension of the Arakawa finite difference C grid, and a dynamical core based on it is currently in development at the UK Met Office \cite{ford2013gung}, due to replace the current finite difference latitude longitude grid discretisation.\\ \\
An important aspect of discretisations in NWP, particularly for climate simulations, is conservation of quantities such as mass and energy. While the former can be conserved using a suitable discretisation of the continuity equation, the latter requires a careful discretisation of all prognostic equations, ensuring that the energy losses and gains are balanced between the discretised terms. A failure to do so may lead to unbalanced transfers between kinetic, potential, and internal energy, which in turn may lead to energy biases in climate models \cite{lucarini2011energetics}. Two ways of guiding this process are given by a variational approach (see e.g. \cite{bauer2019variational, natale2016variational}) and a Hamiltonian framework \cite{morrison1998hamiltonian}. In the latter approach, the Hamiltonian is given by the system's total energy, and the equations are inferred by a Poisson bracket; conservation of energy then follows directly from the bracket's antisymmetry, and any space discretisation maintaining this property will then also conserve energy. For the compressible Euler equations, which form the basic equation set of dynamical cores in NWP, such Hamiltonian based formulations already exist e.g. for hexagonal C-grids \cite{gassmann2013global}, the compatible finite element method \cite{lee2020mixed}, and for general classes of mimetic discretisations \cite{taylor2020energy}.\\ \\
In particular, this framework facilitates the construction of energy conserving higher order numerical methods for transport terms. First, the corresponding terms in the bracket are modified to include the desired numerical method, then the remaining bracket terms are adjusted to retain antisymmetry and therefore energy conservation. For finite element discretisations, the choice of method depends on the underlying space, and for the compatible finite element framework, this requires different formulations for the range of different spaces in use, including continuous and discontinuous ones. Within the Poisson bracket framework, classical DG upwinding and variations thereof for the depth field in the rotating shallow water equations as well as for the thermal field in the thermal shallow water equations have already been introduced in \cite{wimmer2020energy} and \cite{ELDRED20191}, respectively. Further, an energy conserving upwind stabilised discretisetion for the velocity transport term was presented in the context of Lie derivatives for the incompressible Euler equations in \cite{natale2016variational}, which was extended to the shallow water case in \cite{natale2016compatible} and the shallow water case from a Poisson bracket point of view in \cite{wimmer2020energy}. Finally, in the context of continuous or partially continuous finite element spaces, transport terms can be discretised using the Streamline Upwind Petrov Galerkin (SUPG) method \cite{brooks1982streamline}, where a diffusive term is added along the direction of the flow to all test functions. Within the Poisson bracket framework, an energy conserving SUPG method was considered for the evolution of potential vorticity in the shallow water equations in \cite{BAUER2018171}, exploiting the fact that the Poisson bracket term corresponding to the curl part of the velocity transport term's vector-invariant form is antisymmetric in itself. Another Petrov Galerkin type method exploiting the latter fact is given in \cite{lee2021petrov}, where the potential vorticity is diagnosed using Lagrangian trial functions, which are evaluated at downstream locations. \\ \\
These energy conserving upwind-stabilised methods to discretise transport terms can in principle -- with the exception of the SUPG method for potential vorticity -- be implemented readily for the three-dimensional dry compressible Euler equations. However, such a method is still missing for the thermal field transport equation if a continuous (or partially continuous) finite element space is used for the latter field. In particular, this is the case for the finite element equivalent of a Charney-Phillips finite difference grid, where the space's node locations are set to coincide with the velocity space nodes corresponding to the vertical velocity \cite{natale2016compatible}. This choice of thermal field space was recommended in a dispersion property study in \cite{melvin2018choice} for mixed finite element methods in NWP when compared to a fully discontinuous or fully continuous finite element space, and will be used in the UK Met Office's next fluid dynamics component. \\ \\
In this paper, we consider a Poisson bracket together with two Hamiltonians that lead to the thermal rotating shallow water as well as dry compressible Euler equations. We present an almost Poisson bracket based on it that is discretised using the compatible finite element method and includes SUPG for the thermal field transport equation. The resulting space discretisation's energy conservation property as well as the thermal field's qualitative development are verified in numerical test cases, using an energy conserving time discretisation as described in \cite{cohen2011linear}.  For simplicity, while the formulation is valid in three dimensions, we consider two-dimensional scenarios in this paper. This includes vertical slice test cases of the Euler equations with a Charney-Phillips type thermal field space, as well as spherical test cases of the thermal shallow water equations with a fully continuous thermal field. In view of the newly introduced space discretisation's applicability, we also consider the latter in combination with existing energy conserving, upwind-stabilised methods for the density and velocity transport terms, as well as a slightly simplified version of the energy conserving time discretisation. The resulting scheme only contains relatively minor modifications when compared to a corresponding fully upwind-stabilised scheme not derived in an energy conserving context. Next to the aforementioned numerical validations for the energy conserving SUPG method, we will also consider comparisons between the latter two approximately energy conserving and non-energy conserving schemes, which can be run at a comparable computational cost.\\ \\
The rest of the paper is structured as follows: In Section \ref{section_formulation}, we review the Poisson bracket as well as the SUPG method, and present the incorporation of the latter method into the bracket. In Section \ref{section_time_discretisation}, we describe the energy conserving time discretisation and introduce the approximately energy conserving, fully upwind-stabilised scheme. In Section \ref{section_numerical_results}, we present numerical results. Finally, in Section \ref{section_conclusion}, we review the formulation and the corresponding numerical results, and discuss ongoing work.
\section{Formulation} \label{section_formulation}
In this section, we first review the Poisson bracket, which we present in its continuous form. We then derive the resulting sets of equations, considering Hamiltonians that lead to the dry compressible Euler and the thermal shallow water equations. Next, we describe the compatible finite element method to be used to discretise the bracket, including the SUPG method, as well as a mixed SUPG/DG upwind method for the Charney-Phillips type finite element space, as given in \cite{natale2016compatible}. Finally, we incorporate these methods into the Poisson bracket framework, and further review how to include existing energy conserving upwind-stabilisation methods for the velocity and density field related transport terms in the bracket.
\subsection{Hamiltonian framework} \label{section_Poisson_bracket}
Many fluid dynamical equations can be formulated within a Hamiltonian framework, using the system's Hamiltonian $H$, i.e. the total amount of energy, and a Poisson bracket $\{.,.\}$, which is an antisymmetric bilinear form that satisfies the Jacobi identity (for more details, see \cite{shepherd1990symmetries}). The time evolution of any functional $F$ of the prognostic variables is then given by
\begin{equation}
\frac{dF}{dt} = \{F, H\}. \label{Poisson_system}
\end{equation}
Here, we consider as prognostic variables the flow velocity $\mathbf{u}$, density $\rho$, and a thermal type field $\theta$ in a domain $\Omega$. Note that while this choice of variables is common in fluid dynamics, it is non-canonical from a Hamiltonian framework point of view. In particular, the Poisson bracket appearing in \eqref{Poisson_system} will be in a non-canonical form, and here we use a bracket as first introduced in \cite{morrison1980noncanonical}. It is given by
\begin{align}
\begin{split}
\{F,  H\} \coloneqq  - \left\langle \dd{F}{\mathbf{u}}, \mathbf{q} \times \dd{H}{\mathbf{u}} \right\rangle  &+ \left\langle \dd{H}{\rho}, \nabla \cdot \dd{F}{\mathbf{u}} \right\rangle + \left\langle \frac{1}{\rho} \dd{H}{\theta} \nabla \theta, \dd{F}{\mathbf{u}} \right\rangle \label{cts_Euler_bracket}\\
&- \left\langle \dd{F}{\rho}, \nabla \cdot \dd{H}{\mathbf{u}} \right\rangle - \left\langle \frac{1}{\rho} \dd{F}{\theta} \nabla \theta, \dd{H}{\mathbf{u}} \right\rangle,
\end{split}
\end{align}
where $\left\langle . , .\right\rangle$ denotes the $L^2$-inner product over $\Omega$. $\mathbf{q}$ denotes a vorticity type variable, given by
\begin{equation}
\mathbf{q} = (\nabla \times \mathbf{u} + 2\mathbf{\Omega})/\rho, \label{3D_vorticity}
\end{equation}
for rotation vector $\mathbf{\Omega}$. Note that in the two-dimensional test cases considered in the numerical results section below, we have $\mathbf{q} = q\hat{\mathbf{z}}$, for the 2D domain's outward unit vector $\hat{\mathbf{z}}$ and scalar vorticity
\begin{equation}
q = (\nabla^\perp \cdot \mathbf{u} + f)/\rho, \label{potential_vorticity}
\end{equation}
where $\nabla^\perp$ denotes the 2D curl operator $\hat{\mathbf{z}} \times \nabla$, and $f$ is given by the Coriolis force parameter depending on the test case. Finally, given $\mathbf{u} \in V_u(\Omega)$ -- for a suitable space $V_u(\Omega)$ to be defined -- the corresponding functional derivative $dF(\mathbf{u}; \cdot)$ of a given functional $F$ is associated with an element of $V_u(\Omega)$, which is defined weakly by 
\begin{align}
\dd{F}{\mathbf{u}} \in V_u(\Omega) \colon \hspace{1cm} \left\langle \dd{F}{\mathbf{u}} , \mathbf{w} \right\rangle \coloneqq \lim_{\epsilon \rightarrow 0} \frac{1}{\epsilon} \big(F(\mathbf{u} + \epsilon \mathbf{w}, \rho, \theta) -F(\mathbf{u}, \rho, \theta) \big)  && \forall \mathbf{w} \in V_u(\Omega), \label{var_derivative_u}
\end{align}
and similarly for $\rho$ and $\theta$. Note that this definition corresponds to the functional derivative's usual definition,  in the sense that $dF(\mathbf{u}; \cdot)$ -- which is an element of the dual space $V_u(\Omega)'$ -- can be expressed as the 1-form density $\dd{F}{\mathbf{u}} \cdot dx \otimes dV$. Analogous correspondences also hold true for the expressions $\dd{F}{\rho}$ and $\dd{F}{\theta}$ in the density and thermal field spaces $V_\rho(\Omega)$ and $V_\theta(\Omega)$, respectively.\\ \\
In the continuous setting, we consider continuously differentiable spaces, which are then discretised using suitable finite element spaces. The resulting fluid dynamical equations follow from their respective Hamiltonians and functionals $F$ corresponding to weak forms of the prognostic variables. Specifically, we set $F = \left\langle \mathbf{u}, \mathbf{v} \right\rangle$, $F = \left\langle \rho, \psi \right\rangle$, and $F = \left\langle \theta, \chi \right\rangle$, respectively, for arbitrary test functions $\mathbf{v} \in V_u(\Omega)$, $\psi \in V_\rho(\Omega)$, and $\chi \in V_\theta(\Omega)$. For $F = \left\langle \mathbf{u}, \mathbf{v} \right\rangle$, we have
\begin{equation}
\dd{F}{\mathbf{u}} = \mathbf{v}, \; \; \; \dd{F}{\rho} = 0, \; \; \; \dd{F}{\theta} = 0, 
\end{equation} 
so that the Poisson bracket reduces to
\begin{equation}
\{F, H\} = - \left\langle \mathbf{v}, \mathbf{q} \times \dd{H}{\mathbf{u}} \right\rangle  + \left\langle \dd{H}{\rho}, \nabla \cdot \mathbf{v} \right\rangle + \left\langle \frac{1}{\rho} \dd{H}{\theta} \nabla \theta, \mathbf{v} \right\rangle,
\end{equation}
for any test function $\mathbf{v} \in V_u(\Omega)$. Using the Poisson system \eqref{Poisson_system}, we then obtain a momentum equation given by
\begin{align}
\left\langle \mathbf{v}, \mathbf{u}_t \right\rangle  = \frac{dF}{dt} = - \left\langle \mathbf{v}, \mathbf{q}
 \times \dd{H}{\mathbf{u}} \right\rangle  + \left\langle \dd{H}{\rho}, \nabla \cdot \mathbf{v} \right\rangle + \left\langle \frac{1}{\rho} \dd{H}{\theta} \nabla \theta, \mathbf{v} \right\rangle && \forall \mathbf{v} \in V_u(\Omega), \label{weak_Fu_eqn}
\end{align}
where the subscript in $t$ denotes differentiation with respect to time $t$. Similarly, for $\rho$ and $\theta$, we obtain weak versions of the continuity and transport equations of the form
\begingroup
\addtolength{\jot}{2mm}
\begin{align}
&\left\langle \psi, \rho_t \right\rangle = - \left\langle \psi, \nabla \cdot \dd{H}{\mathbf{u}} \right\rangle & \forall \rho \in V_\rho(\Omega), \label{weak_Frho_eqn} \\
&\left\langle \chi, \theta_t \right\rangle = - \left\langle \chi, \frac{1}{\rho} \dd{H}{\mathbf{u}} \cdot \nabla \theta \right\rangle & \forall \chi \in V_\theta(\Omega). \label{weak_Ftheta_eqn}
\end{align}
\endgroup
Note that next to these prognostic equations, the diagnostic vorticity \eqref{potential_vorticity} can also be formulated weakly according to
\begin{align}
\left\langle \zeta, q \rho \right\rangle = - \left\langle \nabla^\perp \zeta, \mathbf{u} \right\rangle + \left\langle\!\!\left\langle \eta, \mathbf{n}^\perp \cdot \mathbf{u} \right\rangle\!\!\right\rangle + \langle \zeta, f\rangle && \forall \zeta \in V_q(\Omega), \label{potential_vorticity_weak}
\end{align}
for $L^2$ inner product $\left\langle\!\left\langle .,. \right\rangle\!\right\rangle$ on the domain's boundary $\partial \Omega$, with outward normal unit vector $\mathbf{n}$, and perpendicular $\mathbf{n}^\perp = \hat{\mathbf{z}} \times \mathbf{n}$. An analogous formulation holds for the 3D version \eqref{3D_vorticity}.\\ \\
Given a system of equations of the form \eqref{weak_Fu_eqn} - \eqref{weak_Ftheta_eqn}, conservation of energy then follows directly from the bracket framework. First, the Poisson system \eqref{Poisson_system} with bracket \eqref{cts_Euler_bracket} can be recovered from \eqref{weak_Fu_eqn} - \eqref{weak_Ftheta_eqn} for any functional $F(\mathbf{u}, \rho, \theta)$ via\footnote{To see this for our choice of notation for the variational derivatives, we note that strictly speaking, the chain rule evaluates to $\frac{dF}{dt} = dF(\mathbf{u}; \mathbf{u}_t) + ... = \left(\dd{F}{\mathbf{u}} \cdot dx \otimes dV\right)(\mathbf{u}_t) + ... = \left\langle \dd{F}{\mathbf{u}}, \mathbf{u}_t \right\rangle + \left\langle \dd{F}{\rho}, \rho_t \right\rangle + \left\langle \dd{F}{\theta}, \theta_t \right\rangle$.}
\begin{equation}
\frac{dF}{dt} = \left\langle \dd{F}{\mathbf{u}}, \mathbf{u}_t \right\rangle + \left\langle \dd{F}{\rho}, \rho_t \right\rangle + \left\langle \dd{F}{\theta}, \theta_t \right\rangle = \{F, H\}. \label{find_dFdt}
\end{equation}
In particular, this also holds for $F=H$, and noting the bracket's antisymmetry, we arrive at
\begin{equation}
\frac{dH}{dt} = \{H, H\} = - \{H, H\} = 0.
\end{equation}
\subsubsection{Euler equations}
A Poisson bracket based formulation of the dry compressible Euler equations can be found e.g. in \cite{gassmann2013global, lee2020mixed}. Note that the Poisson bracket and the skew symmetric operator presented in these papers correspond to a different Poisson bracket to the one considered here, relying on the use of a different set of underlying fields (with a mass weighted thermal field $\Theta$). In our case, the Hamiltonian is given by
\begin{equation}
H(\mathbf{u}, \rho, \theta) = \int_\Omega \big(\frac{\rho}{2} |\mathbf{u}|^2 + g \rho z + c_v \rho \theta \pi \big) \; dx, \label{Euler_Hamiltonian}
\end{equation}\\
for wind velocity $\mathbf{u}$, air density $\rho$, and potential temperature $\theta$. The latter can be seen as a pressure-normalised temperature and is frequently used in numerical weather prediction. For our purposes, we note that unlike the temperature $T$, the potential temperature admits a straightforward transport equation given by \eqref{Euler_th_eqn} below. Note that \eqref{Euler_th_eqn} holds up to source terms and moisture related modifications, which are not considered here. Further, $z$, $g$ and $c_v$ denote the vertical coordinate, gravitational acceleration and specific heat of air at constant volume, respectively. Additionally, $\pi$ denotes the Exner pressure, which is typically considered in conjunction with the potential temperature and can be seen as a normalised pressure. It is given by the ideal gas law as
\begin{equation}
\pi^{\frac{1-\kappa}{\kappa}} = \frac{R}{p_0}\rho \theta, \label{ideal_gas_law}
\end{equation}
for reference pressure $p_0$, ideal gas constant $R$, and non-dimensional parameter $\kappa = R/c_p$, where $c_p = R + c_v$ denotes the specific heat at constant pressure. Finally, in numerical weather prediction, the domain $\Omega$ would typically resemble a spherical shell or a subregion thereof. In the test cases below, we will for simplicity consider (horizontally periodic) rectangular domains, which correspond to local, two-dimensional vertical cross-sections of the aforementioned full domain. In view of the Poisson system \eqref{Poisson_system}, in the non-discretised case the variational derivative of $H$ with respect to $\mathbf{u}$ is given by
\begingroup
\addtolength{\jot}{2mm}
\begin{align}
&\left\langle \dd{H}{\mathbf{u}} , \mathbf{w} \right\rangle = \lim_{\epsilon \rightarrow 0} \frac{1}{\epsilon} \big(H(\mathbf{u} + \epsilon \mathbf{w}, \rho, \theta) - H(\mathbf{u}, \rho, \theta) \big)\\
&\hspace{18mm} = \lim_{\epsilon \rightarrow 0} \frac{1}{\epsilon} \int_\Omega \epsilon \rho \mathbf{u} \cdot \mathbf{w} \; dx = \langle \rho \mathbf{u}, \mathbf{w} \rangle & \forall \mathbf{w} \in V_u(\Omega),\\
\implies &\;\;\; \dd{H}{\mathbf{u}} = \rho \mathbf{u},
\end{align}
\endgroup
and in a similar fashion, we find
\begin{equation}
\dd{H}{\rho} = \frac{1}{2} |\mathbf{u}|^2 + g z +c_p \theta \pi, \hspace{2cm} \dd{H}{\theta} = c_p \rho \pi,
\end{equation}
where we used \eqref{ideal_gas_law} to derive the pressure related expressions. The usual form of the Euler equations then follows from \eqref{weak_Fu_eqn} - \eqref{weak_Ftheta_eqn} and the diagnostic vorticity equation's strong form. For the velocity equation, we have
\begingroup
\addtolength{\jot}{2mm}
\begin{align}
\hspace{-2mm}\left\langle \mathbf{v}, \mathbf{u}_t \right\rangle &= - \left\langle \mathbf{v}, \mathbf{q} \times \dd{H}{\mathbf{u}} \right\rangle  + \left\langle \dd{H}{\rho}, \nabla \cdot \mathbf{v} \right\rangle + \left\langle \frac{1}{\rho} \dd{H}{\theta} \nabla \theta, \mathbf{v} \right\rangle \\
&= - \left\langle \mathbf{v}, \mathbf{q} \times (\rho \mathbf{u}) \right\rangle  + \left\langle \frac{1}{2} |\mathbf{u}|^2 + g z +c_p \theta \pi , \nabla \cdot \mathbf{v} \right\rangle + \left\langle c_p \pi \nabla \theta, \mathbf{v} \right\rangle \\
&= - \left\langle \mathbf{v}, (\nabla \times \mathbf{u} + 2\mathbf{\Omega}) \times \mathbf{u} \right\rangle - \left\langle \frac{1}{2} \nabla |\mathbf{u}|^2 + g \mathbf{k} + c_p \nabla(\theta \pi), \mathbf{v} \right\rangle + \left\langle c_p \pi \nabla \theta, \mathbf{v} \right\rangle\\
&= - \left\langle \mathbf{v}, (\mathbf{u} \cdot \nabla) \mathbf{u} + 2\mathbf{\Omega} \times \mathbf{u} + g \mathbf{k} + c_p \theta \nabla \pi \right\rangle,
\end{align}
\endgroup
for vertical unit vector $\mathbf{k}$ (not to be confused with the domain's outward unit vector $\hat{\mathbf{z}}$), and where we have used integration by parts (assuming $\mathbf{v} \cdot \mathbf{n} = 0$ on $\partial \Omega$) and the vector-invariant identity
\begin{equation}
(\mathbf{u} \cdot \nabla) \mathbf{u} = (\nabla \times \mathbf{u}) \times \mathbf{u} + \frac{1}{2} \nabla |\mathbf{u}|^2. \label{vector_invariant}
\end{equation}
Since this holds for any test function $\mathbf{v}$ in $C^1(\Omega)$, in the case of sufficiently smooth solutions we then arrive at the usual strong form of the momentum equation, given by
\begin{equation}
\mathbf{u}_t + (\mathbf{u} \cdot \nabla) \mathbf{u} + 2\mathbf{\Omega} \times \mathbf{u} + g \mathbf{k} + c_p \theta \nabla \pi = 0, \label{Euler_u_eqn}
\end{equation}
where $2\mathbf{\Omega} \times \mathbf{u}$ corresponds to the Coriolis force, and $c_p \theta \nabla \pi$ is the pressure gradient term expressed in terms of $\theta$ and $\pi$. Similarly, we obtain the continuity and potential temperature transport equations of the form
\begin{align}
&\rho_t + \nabla \cdot (\rho \mathbf{u}) = 0, \\
&\theta_t + \mathbf{u} \cdot \nabla \theta = 0, \label{Euler_th_eqn}
\end{align}
where we used $\dd{H}{\mathbf{u}} = \rho \mathbf{u}$ in \eqref{weak_Frho_eqn} and \eqref{weak_Ftheta_eqn}.
\subsubsection{Thermal shallow water equations}
The thermal shallow water equations have originally been described in \cite{ripa1993conservation}, and here we follow the formulation as given in \cite{ELDRED20191}. The Hamiltonian is given by
\begin{equation}
H(\mathbf{u}, \rho, \theta) = \int_\Omega \left(\frac{\rho}{2}|\mathbf{u}|^2 + \rho\theta\left(\frac{\rho}{2} + b\right)\right)dx, \label{thermal_SWE_Hamiltonian}
\end{equation}
where in this context $\rho$ corresponds to the fluid depth (and for consistency of notation, we keep $\rho$ instead of using the more usual $D$ or $h$). Further, $\theta$ corresponds to the buoyancy, which for our purposes can be seen as a thermal field corresponding to a rescaled temperature (for more details, see e.g. \cite{zeitlin2018geophysical}). In the absence of sources, the buoyancy is governed by the same transport equation as the one used for the potential temperature above. Finally, $b$ denotes the topographic height, and the domain $\Omega$ is given by a sphere. The Hamiltonian variations are now given by
\begin{align}
&\dd{H}{\mathbf{u}} = \rho \mathbf{u}, \hspace{3cm} \dd{H}{\rho} = \frac{1}{2} |\mathbf{u}|^2 + \theta(\rho + b),\\
&\dd{H}{\theta} = \rho\left(\frac{\rho}{2} + b\right),
\end{align}
and we obtain a weak momentum equation of the form
\begingroup
\addtolength{\jot}{2mm}
\begin{align}
\left\langle \mathbf{v}, \mathbf{u}_t \right\rangle =& - \left\langle \dd{F}{\mathbf{u}}, \mathbf{q} \times \dd{H}{\mathbf{u}} \right\rangle  + \left\langle \dd{H}{\rho}, \nabla \cdot \dd{F}{\mathbf{u}} \right\rangle + \left\langle \frac{1}{\rho} \dd{H}{\theta} \nabla \theta, \dd{F}{\mathbf{u}} \right\rangle \\
=& - \left\langle \mathbf{v}, \mathbf{q} \times (\rho \mathbf{u}) \right\rangle  + \left\langle \frac{1}{2} |\mathbf{u}|^2 + \theta(\rho + b), \nabla \cdot \mathbf{v} \right\rangle + \left\langle (\frac{\rho}{2} + b) \nabla \theta, \mathbf{v} \right\rangle \\
=& - \left\langle \mathbf{v}, (\nabla \times \mathbf{u} + 2\mathbf{\Omega}) \times \mathbf{u} \right\rangle - \left\langle \frac{1}{2} \nabla |\mathbf{u}|^2 + \nabla(\theta\big(\rho + b)\big), \mathbf{v} \right\rangle + \left\langle (\frac{\rho}{2} + b) \nabla \theta, \mathbf{v} \right\rangle\\
=& - \left\langle \mathbf{v}, (\mathbf{u} \cdot \nabla) \mathbf{u} +  2 \mathbf{\Omega} \times \mathbf{u} + \theta \nabla (\rho + b) + \frac{\rho}{2} \nabla \theta \right\rangle.
\end{align}
\endgroup
In the non-discretised case, this leads to
\begin{equation}
\mathbf{u}_t + (\mathbf{u} \cdot \nabla) \mathbf{u} + f \mathbf{u}^\perp + \theta \nabla (\rho + b) + \frac{\rho}{2} \nabla \theta = 0,
\end{equation}
noting that given the spherical context of the thermal shallow water equations, we rewrote the Coriolis term as $f\mathbf{u}^\perp$, for Coriolis parameter $f = 2\Omega z/a$ and perpendicular $\mathbf{u}^\perp = \hat{\mathbf{z}} \times \mathbf{u}$, where $a$, $\Omega$, $z$, and $\hat{\mathbf{z}}$ correspond to the sphere's radius, rotation rate, rotational axis coordinate, and outward unit vertical vector, respectively.  Similarly, as in the Euler case, equations for $\rho$ and $\theta$ follow from \eqref{weak_Frho_eqn} and \eqref{weak_Ftheta_eqn}, and since the variational derivative $\dd{H}{\mathbf{u}}$ is the same here, the resulting equations are given as before by
\begin{align}
&\rho_t + \nabla \cdot (\rho \mathbf{u}) = 0, \\
&\theta_t + \mathbf{u} \cdot \nabla \theta = 0.
\end{align}
\subsection{Space discretisation}
In this section, we describe the space discretisation for the Poisson bracket introduced above. First, we discuss our choice of compatible finite element spaces for the fields $\mathbf{u}$, $\rho$, and $q$, as well as the additional space for the thermal field $\theta$. Next, we describe the SUPG method used for the stabilisation of the thermal field transport equation and present its incorporation into the Poisson bracket framework, which forms the core part of this paper. Finally, we discuss how to include upwind-stabilisation for the transport terms corresponding to $\mathbf{u}$ and $\rho$ in this bracket, as presented in \cite{wimmer2020energy}. 
\subsubsection{Choice of compatible finite element spaces} \label{section_compatible_FEM}
\begin{figure}[ht]
\begin{center}
\includegraphics[width=0.7\textwidth]{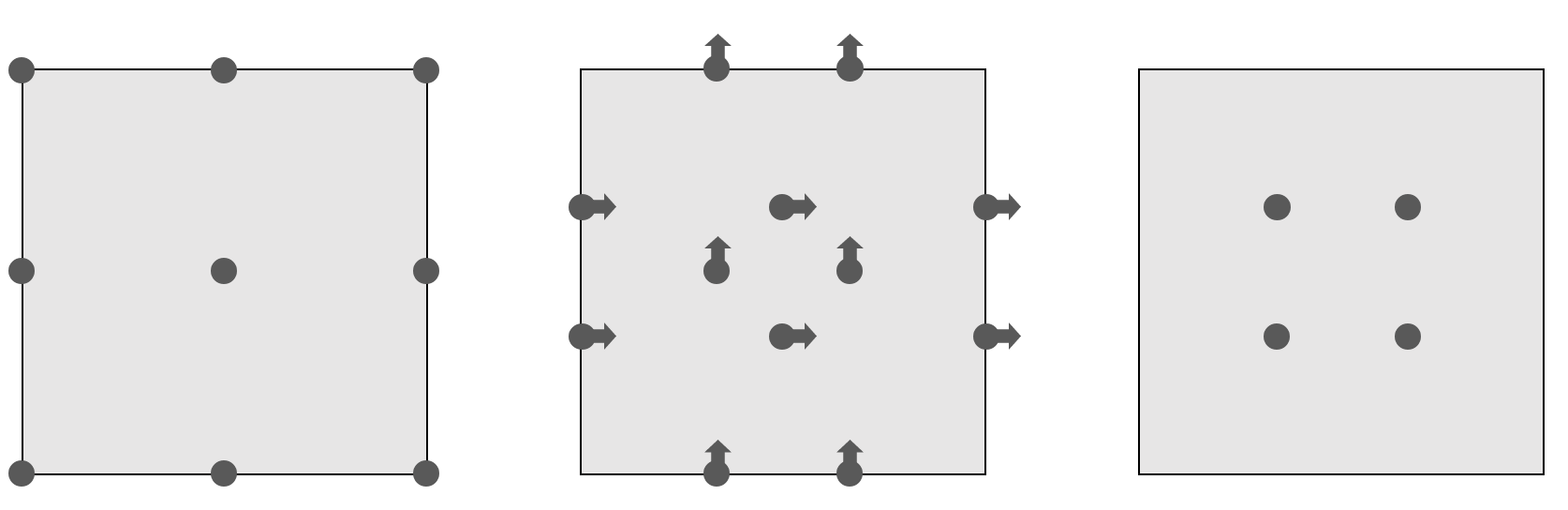}
\caption{Vertical slice reference elements for next to lowest order compatible spaces for $\mathbb{V}_q$, $\mathbb{V}_u$, and $\mathbb{V}_\rho$. Dots and arrows denote degrees of freedom and degrees of freedom multiplied by horizontal/vertical unit vectors, respectively.}\label{Vertical_elements}
\end{center}
\end{figure}

\begin{figure}[ht]
\begin{center}
\includegraphics[width=0.7\textwidth]{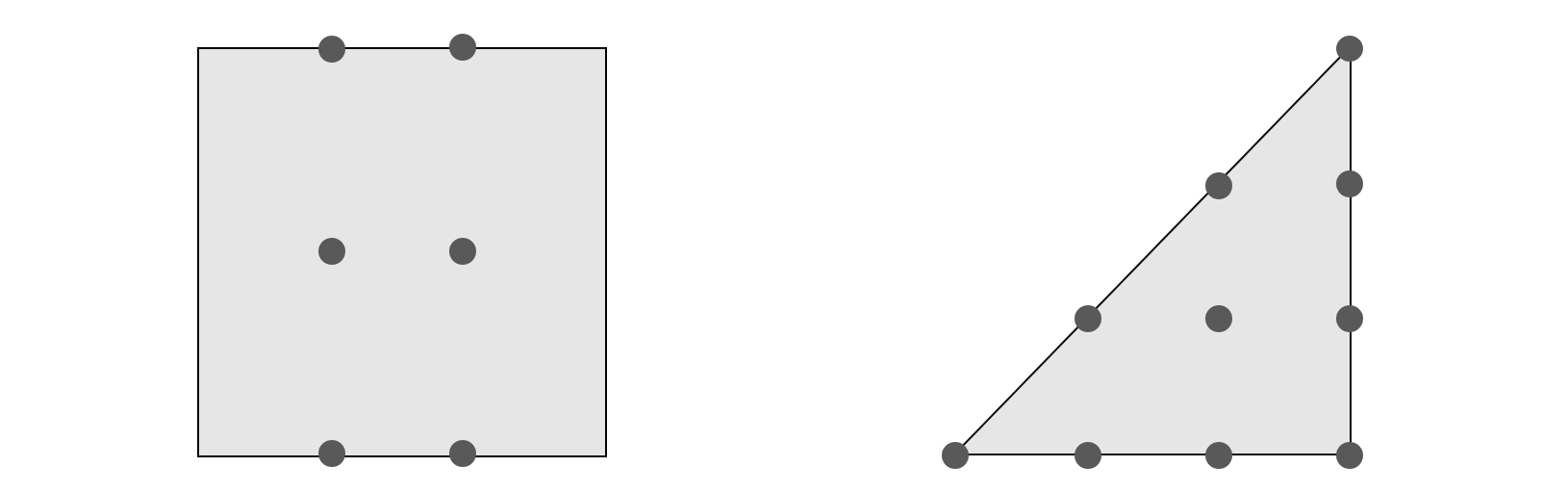}
\caption{Reference elements for next to lowest order Charney-Phillips type finite element space and $P_3$. Dots denote degrees of freedom. The former is used for the Euler equations in a vertical slice of the atmosphere; the latter corresponds to a triangular mesh of the sphere used for the thermal shallow water equations.}\label{Temperature_elements}
\end{center}
\end{figure}

To discretise \eqref{weak_Fu_eqn} - \eqref{weak_Ftheta_eqn} and the weak diagnostic vorticity equation \eqref{potential_vorticity_weak}, we first note that suitable spaces for the vorticity, velocity, and density before discretisation are given by
\begin{equation}
V_q(\Omega) = H^1(\Omega), \hspace{1cm} V_u(\Omega) = H(\text{div}; \Omega), \hspace{1cm} V_\rho(\Omega) = L^2(\Omega),
\end{equation}
since in this case all differential operations in the equations are well-defined (up to $\nabla \theta$). These spaces form a so-called Hilbert complex
\begin{displaymath}
    \xymatrix{
        H^1(\Omega) \ar[r]^{\nabla^\perp} & H(\text{{\normalfont div}};\Omega) \ar[r]^{\nabla \cdot} & L^2(\Omega)},
\end{displaymath}
that is the differential operators occurring in the previously described equations map one space onto another. This motivates the compatible finite element framework, where we consider spaces for $q$, $\mathbf{u}$, and $\rho$, given by $\mathbb{V}_q$, $\mathbb{V}_u$, and $\mathbb{V}_\rho$, respectively, such that
\begin{displaymath}
    \xymatrix{
        V_q(\Omega) \ar[r]^{\nabla^\perp} \ar[d]^{\pi^0} & V_u(\Omega) \ar[d]^{\pi^1} \ar[r]^{\nabla \cdot} & V_\rho(\Omega) \ar[d]^{\pi^2} \\
        \mathbb{V}_q(\Omega) \ar[r]^{\nabla^\perp}       & \mathbb{V}_u(\Omega) \ar[r]^{\nabla \cdot}  & \mathbb{V}_\rho(\Omega),}\label{DeRham}
\end{displaymath}\\
where $\pi^0$, $\pi^1$, and $\pi^2$ are given by suitable bounded projections such that the diagram commutes. This choice of finite element spaces can be shown to ensure a number of discretisation properties relevant to numerical weather prediction, including stability of the pressure gradient operator, as well as exact equilibria\footnote{The equilibria are given by steady geostrophic states, that is steady state solutions in which the Coriolis and pressure gradient terms are balanced.} for the linearised rotating shallow water equations \cite{natale2016compatible}.\\ \\
In the case of the Euler equations, where we consider a (horizontally periodic) vertical slice of the atmosphere, we use regular quadrilateral meshes and set the velocity space equal to the Raviart-Thomas space $RT_k$, for biquadratic ($k=2)$ and third order ($k=3$) finite elements. The compatible density and vorticity spaces are then given by the continuous Galerkin space $Q_k$, as well as the discontinuous Galerkin space $Q_k^{DG}$, respectively. The corresponding compatible finite elements for $k=2$ are depicted in Figure \ref{Vertical_elements}. Note that in the implementation discussed in the numerical results section below, the meshes are given by vertically extruded ones, and the underlying finite elements are implemented as tensor products of one-dimensional finite elements. This readily extends to three-dimensional domains, where we would use two-dimensional compatible base elements in the tensor product (for further details, see \cite{natale2016compatible}). The discretisation described below is still valid in the 3D case, except that the diagnostic scalar vorticity equation would be replaced by its corresponding vector field form.\\
\\
For the thermal shallow water equations, where we consider triangular spherical meshes, we follow \cite{wimmer2020energy} and use the degree two Brezzi-Douglas-Marini finite element space $BDM_2$ for the velocity field. The compatible depth and vorticity spaces are then given by $P_1^{DG}$ and $P_3$, i.e. the first and third (polynomial) order discontinuous and continuous triangular Galerkin spaces, respectively.\\ \\
Finally, we need to specify finite element spaces for the thermal field $\theta$. In the case of the thermal shallow water equations, we use the continuous space $\mathbb{V}_q$, given by $P_3$. For the Euler equations, we consider a Charney-Phillips type space by choosing the finite element nodes to coincide with the velocity space nodes corresponding to wind in the vertical direction. The resulting elements for the thermal field space in the case of the Euler and thermal shallow water equations are given in Figure \ref{Temperature_elements}.
\\
In terms of the equations derived by the Poisson bracket in Section \ref{section_Poisson_bracket}, we obtain a space discretisation given by
\begingroup
\addtolength{\jot}{2mm}
\begin{align}
&\left\langle \mathbf{w}, \mathbf{u}_t \right\rangle+ \left\langle \mathbf{w}, \mathbf{q} \times \dd{H}{\mathbf{u}} \right\rangle  - \left\langle \dd{H}{\rho}, \nabla \cdot \mathbf{w} \right\rangle-  \left\langle \frac{1}{\rho} \dd{H}{\theta} \nabla \theta, \mathbf{w} \right\rangle = 0 & \forall \mathbf{w} \in \mathbb{V}_u, \label{dscr_weak_Fu_eqn}\\
&\left\langle \phi, \rho_t \right\rangle + \left\langle \phi, \nabla \cdot \dd{H}{\mathbf{u}} \right\rangle = 0 \label{dscr_weak_Frho_eqn} &\forall \phi \in \mathbb{V}_\rho, \\
&\left\langle \gamma, \theta_t \right\rangle + \left\langle \gamma, \frac{1}{\rho} \dd{H}{\mathbf{u}} \cdot \nabla \theta \right\rangle =0 &\forall \gamma \in \mathbb{V}_\theta, \label{dscr_weak_Ftheta_eqn}
\end{align}
\endgroup
together with a discretised diagnostic scalar vorticity equation given by
\begin{align}
\left\langle \eta, q\rho \right\rangle = - \left\langle \nabla^\perp \eta, \mathbf{u} \right\rangle + \left\langle\!\!\left\langle \eta, \mathbf{n}^\perp \cdot \mathbf{u} \right\rangle\!\!\right\rangle + \left\langle \eta, f \right\rangle && \forall \eta \in \mathbb{V}_q. \label{discr_diag_q}
\end{align}
Note that we further consider a discrete Hamiltonian, in the sense that it is now a functional of the finite element fields $\mathbf{u} \in \mathbb{V}_u$, $\rho \in \mathbb{V}_\rho$, and $\theta \in \mathbb{V}_\theta$. In particular, the variational derivative \eqref{var_derivative_u} is then defined with respect to the discrete space $\mathbb{V}_u$, and similarly for the variations in $\rho$ and $\theta$. The variations evaluate as projections, and e.g. for the Hamiltonian corresponding to the Euler equations, we find
\begingroup
\addtolength{\jot}{2mm}
\begin{align}
&\left\langle \dd{H}{\mathbf{u}}, \mathbf{w} \right\rangle = \left\langle \rho\mathbf{u}, \mathbf{w} \right\rangle &\forall \mathbf{w} \in \mathbb{V}_u,\\
\implies &\;\;\;\; \dd{H}{\mathbf{u}} = P_{\mathbb{V}_u}(\rho \mathbf{u}),
\end{align}
\endgroup
and similarly
\begin{equation}
\dd{H}{\rho} = P_{\mathbb{V}_\rho}\left(\frac{1}{2} |\mathbf{u}|^2 + g z +c_p \theta \pi\right), \hspace{2cm} \dd{H}{\theta} = P_{\mathbb{V}_\theta}(c_p \rho \pi),
\end{equation}
where $P_\mathbb{V}$ denotes the $L^2$-projection into $\mathbb{V}$. Further, we note that the Poisson bracket structure still holds in the discretised case, i.e. the discretised evolution equations \eqref{dscr_weak_Fu_eqn} - \eqref{dscr_weak_Ftheta_eqn} can be used for evaluating $dF/dt$ as in \eqref{find_dFdt}, leading to a Poisson system (as in \eqref{Poisson_system}) with a bracket whose form is equal to the continuous one, given by \eqref{cts_Euler_bracket}. While it is not clear if the discretised bracket still satisfies the Jacobi identity, it is still antisymmetric, implying that energy is conserved. We therefore find that the discretised bracket is at least an almost Poisson bracket, i.e. a bilinear antisymmetric form.
\subsubsection{SUPG stabilisation} \label{section_SUPG}
To review the SUPG method, we start from the thermal field transport equation in its continuous form, given by
\begin{equation}
\theta_t + \mathbf{u} \cdot \nabla \theta = 0. \label{advection_equation}
\end{equation}
In a finite element setting, we would then space discretise using its weak form
\begin{align}
\left\langle \gamma, \theta_t + \mathbf{u} \cdot \nabla \theta \right\rangle = 0 && \forall \gamma \in \mathbb{V}_\theta. \label{weak_advection_equation}
\end{align}
However, for standard continuous Galerkin spaces, the advection term tends to produce spurious oscillations, which can be counteracted by introducing diffusion in the direction of the flow \cite{kuzmin2010guide}. This is achieved by replacing the test functions by
\begin{equation}
\gamma \rightarrow \gamma + \tau \mathbf{u} \cdot \nabla \gamma, \label{gamma_upw}
\end{equation}
for a suitable coefficient $\tau$. We then arrive at
\begin{align}
\left\langle \gamma + \tau \mathbf{u} \cdot \nabla \gamma, \theta_t + \mathbf{u} \cdot \nabla \theta \right\rangle = 0 && \forall \gamma \in \mathbb{V}_\theta, \label{theta_SUPG_eqn}
\end{align}
and this form of upwinding is commonly referred to as the Streamline Upwind Petrov Galerkin (SUPG) method \cite{kuzmin2010guide}, and was first introduced in \cite{brooks1982streamline}. We find that the modified test function of the transport term acts as diffusion along the velocity field $\mathbf{u}$, noting that for $\gamma = \theta \in \mathbb{V}_\theta$, the SUPG equation \eqref{theta_SUPG_eqn} leads to
\begin{align}
\frac{1}{2} \frac{d}{dt} \|\theta\|_2^2 = - \left\langle \theta, \mathbf{u} \cdot \nabla \theta \right\rangle - \| \sqrt{\tau} \mathbf{u} \cdot \nabla \theta \|_2^2 - \left\langle \tau \mathbf{u} \cdot \nabla \theta, \theta_t \right\rangle, \label{theta_dissipation}
\end{align}
where the additional indefinite SUPG term (i.e. the last term) is kept small by a suitable choice of $\tau$, which is often given by $\tau = \Delta x_h/ (2|\mathbf{u}|)$, for local mesh size $\Delta x_h$ \cite{kuzmin2010guide}.  Since all test functions were modified to include an upwind term, the resulting weak equation \eqref{theta_SUPG_eqn} is still clearly consistent, in the sense that an exact solution to the continuous equation \eqref{advection_equation} also solves the weak equation.
\\ \\
The upwind-stabilised transport equation \eqref{theta_SUPG_eqn} can be used for fully continuous spaces $\mathbb{V}_\theta$, as is the case for our choice of thermal field space for the thermal shallow water equations. For spaces that are not fully continuous, such as the Charney-Phillips type finite element space defined above, further care needs to be taken for the discontinuous components. In our case, we follow \cite{natale2016compatible}, and apply the standard DG upwind method \cite{kuzmin2010guide} in the horizontal (i.e. discontinuous) direction, and restrict the SUPG setup to the vertical direction. To this end, we first integrate the weak form \eqref{weak_advection_equation} by parts and arrive at
\begin{align}
\left\langle \gamma, \theta_t \right\rangle - \left\langle \nabla \cdot (\mathbf{u}\gamma), \theta \right\rangle + \int_{\Gamma_v} \left[\!\left[\mathbf{u} \gamma \right]\!\right] \tilde{\theta} \; dS = 0 && \forall \gamma \in \mathbb{V}_\theta, \label{CP_DG_upw_incomplete}
\end{align}
where $\Gamma_v$ denotes the set of all vertical interior facets of the underlying mesh. The jump operator $\left[\!\left[.\right]\!\right]$ (for vectors $\mathbf{v}$ and scalars $\psi$, respectively) and upwind value $\tilde{\theta}$ are defined by
\begin{equation}
\begin{aligned}
&\left[\!\left[\mathbf{v}\right]\!\right] = \mathbf{v}^+ \cdot \mathbf{n}^+ + \mathbf{v}^- \cdot \mathbf{n}^-,\\
&\left[\!\left[\psi\right]\!\right] = \psi^+ - \psi^-,
\end{aligned}
 \hspace{2cm} \tilde{\theta}=
     \begin{cases}
     	\begin{split}
       	\theta^+ \; \; \; \; \text{if } \mathbf{u} \cdot \mathbf{n}^+ < 0,\\
	\theta^- \; \; \; \; \text{otherwise}, \label{def_upwind_theta} \; \; \; \;
	\end{split}
     \end{cases}
\end{equation}
for a given facet's normal unit vector $\mathbf{n}$, and noting that the two sides of each mesh facet are arbitrarily denoted by + and - (and hence $\mathbf{n}^+ = -\mathbf{n}^-$). Before incorporating the upwind modification \eqref{gamma_upw} to the test functions in \eqref{CP_DG_upw_incomplete}, we need to integrate by parts again to avoid applying the differential operator $\nabla$ to the upwinded part $\tau \mathbf{u} \cdot \nabla \gamma$, as the double differentiation may not be well-defined for $\gamma \in \mathbb{V}_\theta$. Further, we restrict the SUPG method to the vertical direction by using modified test functions of the form
\begin{equation}
\gamma \rightarrow \gamma + \tau (\mathbf{k} \cdot {\mathbf{u}})(\mathbf{k} \cdot \nabla \gamma) \eqqcolon \gamma_{(\mathbf{k} \cdot \mathbf{u}) \mathbf{k}}, \label{gamma_upw_CP}
\end{equation}
for vertical unit vector $\mathbf{k}$. The resulting, fully upwinded transport equation is then given by
\begin{align}
\left\langle \gamma_{(\mathbf{k} \cdot \mathbf{u}) \mathbf{k}}, \theta_t \right\rangle + \left\langle \gamma_{(\mathbf{k} \cdot \mathbf{u}) \mathbf{k}}, \mathbf{u} \cdot \nabla \theta \right\rangle + \int_{\Gamma_v} \big( \left[\!\left[\mathbf{u} \gamma_{(\mathbf{k} \cdot \mathbf{u}) \mathbf{k}} \right]\!\right] \tilde{\theta} - \left[\!\left[\mathbf{u} \gamma_{(\mathbf{k} \cdot \mathbf{u}) \mathbf{k}} \theta\right]\!\right] \big)dS = 0 && \forall \gamma \in \mathbb{V}_\theta, \label{theta_SUPG_eqn_CP}
\end{align} \\
noting that the second integration by parts leads to another facet integral term. Further, note that given the partially discontinuous Charney-Phillips type finite element space, the gradient term in \eqref{theta_SUPG_eqn_CP} is considered element-wise. In the next section, we will incorporate the SUPG methods used to stabilise equations \eqref{theta_SUPG_eqn} and \eqref{theta_SUPG_eqn_CP} in the energy conserving bracket \eqref{cts_Euler_bracket}. To simplify notation, we write the modified test functions \eqref{gamma_upw} and \eqref{gamma_upw_CP} as well as the transport terms in \eqref{theta_SUPG_eqn} and \eqref{theta_SUPG_eqn_CP} as single expressions given in the definition below, where CG and CP denote continuous Galerkin and Charney-Phillips, respectively.
\begin{definition} \label{def_gamma_upw_L_operator}
Consider a velocity field element $\mathbf{u} \in \mathbb{V}_u$ and a thermal field space element $\gamma \in \mathbb{V}_\theta$. The SUPG contribution corresponding to $\gamma$ is then defined by
\begin{equation}
S(\mathbf{u}; \gamma) =
    \begin{cases}
     	\begin{split}
       	\mathbf{u} \cdot \nabla \gamma \; \; \; \; &\text{CG type } \mathbb{V}_\theta,\\
	(\mathbf{k} \cdot \mathbf{u})(\mathbf{k} \cdot \nabla \gamma) \; \; \; \; &\text{CP type } \mathbb{V}_\theta. \label{def_upwind_gamma}
	\end{split}
    \end{cases}
\end{equation}
Further, consider another thermal field space element $\theta \in \mathbb{V}_\theta$ and a function $\sigma$, which may correspond to $\gamma$ or its modified version $\gamma + \tau S(\mathbf{u}; \gamma)$. The discrete thermal field transport expression $L(\mathbf{u}, \theta; \sigma)$ is defined by
\begin{equation}
\hspace{-2mm} L(\mathbf{u}, \theta; \sigma) \!=\!
     \begin{cases}
     	\begin{split}
       	&\!\! \left\langle \sigma , \mathbf{u} \cdot \nabla \theta \right\rangle & \text{CG type } \mathbb{V}_\theta,\\
	&\!\! \left\langle \sigma, \mathbf{u} \cdot \nabla \theta \right\rangle \!+\!\! \int_{\Gamma_v} \!\!\! \big( \left[\!\left[\mathbf{u} \sigma \right]\!\right] \tilde{\theta} - \left[\!\left[\mathbf{u} \sigma \theta\right]\!\right] \big)dS & \text{CP type } \mathbb{V}_\theta.\label{th_adv_operator}
	\end{split}
     \end{cases}
\end{equation}
\end{definition}
Using this definition, we may rewrite the discretised thermal field transport equation as
\begin{align}
\left\langle \gamma, \theta_t \right\rangle + L(\mathbf{u}, \theta; \gamma) = 0 && \forall \gamma \in \mathbb{V}_\theta, \label{weak_advection_equation_L}
\end{align}
and its SUPG stabilised versions \eqref{theta_SUPG_eqn} and \eqref{theta_SUPG_eqn_CP} as a single equation of the form
\begin{align}
\left\langle \gamma + \tau S(\mathbf{u}; \gamma), \theta_t \right\rangle + L\big(\mathbf{u}, \theta; \gamma + \tau S(\mathbf{u}; \gamma)\big) = 0 && \forall \gamma \in \mathbb{V}_\theta. \label{theta_SUPG_eqn_cpt}
\end{align}
Finally, we note that the first argument in $L$ corresponds to the advecting velocity, while the first one in $S$ does to the velocity used for upwinding. Depending on the time discretisation, we will find that these may be distinct in the almost Poisson bracket setup.
\subsubsection{SUPG in the discretised bracket} \label{section_SUPG_bracket}
It remains to incorporate the SUPG method, as described in Section \ref{section_SUPG}, into the discretised almost Poisson bracket. Before moving on to relevant definitions for the updated bracket, we outline the course of thought leading to it. Once the bracket has been defined, we show how to derive the corresponding momentum and thermal field equations, and demonstrate that the total energy is still conserved. Finally, we discuss the impact of the term in the bracket that is antisymmetric to the term corresponding to the SUPG method, and further include a remark on transport stabilisation for the other prognostic variables.\\ \\
Without upwinding, the bracket is given by \eqref{cts_Euler_bracket}, such that $q \in \mathbb{V}_q$, $\rho \in \mathbb{V}_\rho$, $\theta \in \mathbb{V}_\theta$, and the Hamiltonian variations are functions of the corresponding finite element spaces. We aim to replace the corresponding thermal field advection equation \eqref{dscr_weak_Ftheta_eqn} by the SUPG version \eqref{theta_SUPG_eqn_cpt},  i.e.
\begin{equation}
\left\langle \gamma, \theta_t \right\rangle = - L\left(\frac{1}{\rho}\dd{H}{\mathbf{u}}, \theta; \gamma \right)  \; \; \; \rightarrow \; \; \; \left\langle \gamma + \tau S(\mathbf{u}; \gamma), \theta_t \right\rangle = - L\left(\frac{1}{\rho}\dd{H}{\mathbf{u}}, \theta; \gamma + \tau S(\mathbf{u}; \gamma) \right). \label{SUPG_theta_eqn_bracket_form}
\end{equation}
First, note that since the advecting velocity in the non-stabilised thermal field transport equation is given by $\dd{H}{\mathbf{u}}/\rho$, we anticipate that this will also be the case when the SUPG method is included. The non-stabilised version (to the left in \eqref{SUPG_theta_eqn_bracket_form}) was derived from choosing $F = \left\langle \gamma, \theta \right\rangle$ in the Poisson system \eqref{Poisson_system}, suggesting that for the stabilised one we could instead pick a functional of the form
\begin{equation}
F = \left\langle \gamma + \tau S(\mathbf{u}; \gamma), \theta \right\rangle.
\end{equation}
However, this is impractical, since then the time derivative of $F$ will include a time derivative in the upwinding velocity $\mathbf{u}$ appearing in the upwind contribution $S$. To avoid this, we introduce an operator corresponding to the SUPG-modified test function $\gamma + \tau S(\mathbf{u}; \gamma)$, as defined below.
\begin{definition} \label{definition_s}
Consider the SUPG parameter $\tau$, a velocity field $\mathbf{u} \in \mathbb{V}_u$, and a thermal field $\gamma \in \mathbb{V}_\theta$. Then the SUPG operator $s \colon \mathbb{V}_\theta \rightarrow \mathbb{V}_\theta$, parameterised by $\tau$ and $\mathbf{u}$, is defined such that
\begin{align}
s = s(\tau, \mathbf{u}; \gamma)\colon \hspace{1cm} \left\langle s + \tau S(\mathbf{u}; s), \sigma \right\rangle = \langle \gamma,  \sigma \rangle && \forall \sigma \in \mathbb{V}_\theta. \label{def_s_form}
\end{align}
\end{definition}
\begin{remark} \label{remark_app_A}
Note that the computation of $s$ corresponds to inverting a mass matrix that is modified to also include an SUPG contribution. The invertibility of such modified mass matrices has been studied e.g. in \cite{bochev2004stability} for divergence-free flow fields and Lagrange elements. The extension of the proofs given therein to the non-divergence-free flow fields and finite elements used here is given in Appendix \ref{appendix_A}.
\end{remark}
Using $s$, we may keep a functional $F$ of the form $F = \langle \gamma, \theta \rangle$, leading to a Poisson system formulation given by
\begin{equation}
\left\langle s + \tau S(\mathbf{u}; s), \theta_t \right\rangle = \langle \gamma, \theta_t \rangle = \frac{dF}{dt} = \{F, H\}. \label{Poisson_system_s}
\end{equation}
The left-hand side of \eqref{Poisson_system_s} then corresponds to the left-hand side of our desired SUPG stabilised thermal field equation \eqref{SUPG_theta_eqn_bracket_form}, where $s$ now plays the role of the test function (rather than $\gamma$). To also obtain the required right-hand side in \eqref{SUPG_theta_eqn_bracket_form}, we require a Poisson bracket that yields, for $F = \langle \gamma, \theta \rangle$, a transport term of the form
\begin{equation}
\{F, H\} = - L\left(\frac{1}{\rho}\dd{H}{\mathbf{u}}, \theta; s + \tau S(\mathbf{u}; s) \right),
\end{equation}
recalling that $s =s(\tau, \mathbf{u}; \gamma)$. To achieve this, we simply replace the antisymmetric pair corresponding to the thermal field terms in the bracket \eqref{cts_Euler_bracket}, motivating a modified bracket as given in Definition \ref{def_SUPG_bracket} below.
\begin{definition} \label{def_SUPG_bracket}
The SUPG almost Poisson bracket is given by 
\begingroup
\addtolength{\jot}{4mm}
\begin{align}
\begin{split}
\hspace{-3mm} \{F, H\} = &- \left\langle \dd{F}{\mathbf{u}}, \mathbf{q} \times \dd{H}{\mathbf{u}} \right\rangle \\
& + \left\langle \dd{H}{\rho}, \nabla \cdot \dd{F}{\mathbf{u}} \right\rangle + L\left(\frac{1}{\rho}\dd{F}{\mathbf{u}}, \theta; s\left(\tau, \mathbf{u}; \dd{H}{\theta}\right) + \tau S\left(\mathbf{u}; s\left(\tau, \mathbf{u}; \dd{H}{\theta}\right)\right)\right) \\
& - \left\langle \dd{F}{\rho}, \nabla \cdot \dd{H}{\mathbf{u}} \right\rangle - L\left(\frac{1}{\rho}\dd{H}{\mathbf{u}}, \theta; s\left(\tau, \mathbf{u}; \dd{F}{\theta}\right) + \tau S\left(\mathbf{u}; s\left(\tau, \mathbf{u}; \dd{F}{\theta}\right)\right)\right), \label{SUPG_bracket}
\end{split}
\end{align}
\endgroup\\
where the expressions $L$ and $S$ are given by Definition \ref{def_gamma_upw_L_operator}, and $s$ is defined according to Definition \ref{definition_s}.
\end{definition}
Note that the above bracket is antisymmetric by construction. Further, it is bilinear since $L$ is bilinear with respect to its first and third arguments (the advecting velocity and modified test function), and $S$ and $s$ are linear in their last argument as well. In particular, like the discretised version of \eqref{cts_Euler_bracket}, the SUPG modified bracket \eqref{SUPG_bracket} is indeed an almost Poisson bracket.\\ \\
Additionally, note that if we set $\tau = 0$ and ignore the modifications to the transport operator $L$ related to the Charney-Phillips grid, then the bracket reduces to (the discrete equivalent of) the original bracket \eqref{cts_Euler_bracket}. This holds true since in this case the $S$ expressions vanish, and the operator $s$ reduces to a projection into the thermal field space, which is equal to the identity operator for $\dd{H}{\theta} \in \mathbb{V}_\theta$ and $\dd{F}{\theta} \in \mathbb{V}_\theta$. \\ \\
In Proposition \ref{prop_1}, we confirm that this setup leads to an SUPG stabilised thermal field equation of the form \eqref{theta_SUPG_eqn_cpt}.
%
%
%
\begin{prop} \label{prop_1}
Consider the space discretised almost Poisson bracket given by Definition \ref{def_SUPG_bracket}. Then the thermal field transport equation based on the Poisson system \eqref{Poisson_system} with this bracket is given by
\begin{align}
\left\langle s + \tau S(\mathbf{u}; s), \theta_t \right\rangle + L\left(\frac{1}{\rho}\dd{H}{\mathbf{u}}, \theta; s + \tau S(\mathbf{u}; s)\right) = 0 && \forall s \in \mathbb{V}_\theta. \label{SUPG_theta_ec}
\end{align}
\textbf{Proof}. 
We consider the Poisson system with the usual functional of the form $F = \left\langle \gamma, \theta \right\rangle$, for any $\gamma \in \mathbb{V}_\theta$. In this case, as before all parts of the bracket except for the one corresponding to thermal field transport (i.e. the last term in \eqref{SUPG_bracket}) are zero, and we are left with
\begin{align}
\left\langle s + S\big(\tau, \mathbf{u}; s\big), \theta_t \right\rangle = \left\langle \gamma, \theta_t \right\rangle = \frac{dF}{dt} = \{F, H\} =   - L\left(\frac{1}{\rho}\dd{H}{\mathbf{u}}, \theta; s + S\big(\tau, \mathbf{u}; s\big)\right) && \forall \gamma \in \mathbb{V}_\theta, \label{theta_eqn_forall_gamma}
\end{align}
where $s = s(\tau, \mathbf{u}; \gamma)$ is given by \eqref{def_s_form}. Note that $s$ spans the thermal field space $\mathbb{V}_\theta$, since it corresponds to an invertible operator (see Appendix \ref{appendix_A}). In particular, since \eqref{theta_eqn_forall_gamma} holds for all $\gamma \in \mathbb{V}_\theta$, it also does for all $s = s(\tau, \mathbf{u}; \gamma) \in \mathbb{V}_\theta$.
\hfill$\square$
\end{prop}
In Proposition \ref{prop_2}, we derive the momentum equation arising from the SUPG bracket \eqref{SUPG_bracket}, which now also contains SUPG related terms due to the bracket's antisymmetric structure.
%
%
%
\begin{prop} \label{prop_2}
Consider the space discretised almost Poisson bracket given by Definition \ref{def_SUPG_bracket}. Then the momentum equation based on the Poisson system \eqref{Poisson_system} with this bracket is given by
\begin{align}
\begin{split}
\left\langle \mathbf{w}, \mathbf{u}_t \right\rangle &+ \left\langle \mathbf{w}, \mathbf{q} \times \dd{H}{\mathbf{u}} \right\rangle - \left\langle \dd{H}{\rho}, \nabla \cdot \mathbf{w} \right\rangle \\
&- L\left(\frac{1}{\rho}\mathbf{w}, \theta; s\left(\tau, \mathbf{u}; T\right) + \tau S\left(\mathbf{u}; s\left(\tau, \mathbf{u}; T\right)\right)\right)  = 0 \hspace{15mm} \forall \mathbf{w} \in \mathbb{V}_u,
\end{split}
\end{align}
where $T$ is given by
\begin{align}
T=
\begin{cases}
\;\;\rho(\frac{\rho}{2} + b) &\text{Thermal shallow water equations,}\\
\;\;c_p \rho \pi &\text{Euler equations.}
\end{cases}
\end{align}
\textbf{Proof.} Considering the Poisson system \eqref{Poisson_system} together with the SUPG bracket \eqref{SUPG_bracket} and a functional of the form $F = \langle \mathbf{w}, \mathbf{u} \rangle$ (for any $\mathbf{w} \in \mathbb{V}_u$), we arrive at a momentum equation of the form
\begin{align}
\begin{split}
\left\langle \mathbf{w}, \mathbf{u}_t \right\rangle &+ \left\langle \mathbf{w}, \mathbf{q} \times \dd{H}{\mathbf{u}} \right\rangle - \left\langle \dd{H}{\rho}, \nabla \cdot \mathbf{w} \right\rangle \\
&- L\left(\frac{1}{\rho}\mathbf{w}, \theta; s\left(\tau, \mathbf{u}; \dd{H}{\theta}\right) + S\left(\tau, \mathbf{u}; s\left(\tau, \mathbf{u}; \dd{H}{\theta}\right)\right)\right)  = 0 \hspace{15mm} \forall \mathbf{w} \in \mathbb{V}_u.
\end{split}
\end{align}
Note that $s$ is derived according to \eqref{def_s_form}, i.e.
\begin{align}
\left\langle s + S\big(\tau, \mathbf{u}; s\big), \sigma \right\rangle = \left\langle \dd{H}{\theta},  \sigma \right\rangle && \forall \sigma \in \mathbb{V}_\theta.
\end{align}
Substituting the expression for $\dd{H}{\theta}$ for the thermal shallow water equations, we find
\begin{align}
\left\langle \dd{H}{\theta}, \sigma \right\rangle = \left\langle \rho(\frac{\rho}{2} + b), \sigma \right\rangle && \forall \sigma \in \mathbb{V}_\theta,
\end{align}
and similarly for the Euler equations. In other words,
\begin{equation}
s\left(\tau, \mathbf{u}; \dd{H}{\theta}\right) = s\left(\tau, \mathbf{u}; T\right),
\end{equation}
as required. \hfill$\square$
\end{prop}
When comparing the momentum equation given in Proposition \ref{prop_2} above with the original discretised bracket's momentum equation \eqref{dscr_weak_Fu_eqn}, we find that an additional computation is required to obtain $s$, while the computation for the projection $\dd{H}{\theta}$ is not required anymore. In particular, since both computations are in the thermal field space $\mathbb{V}_\theta$, the overall computational cost is similar.\\ \\
Using Propositions \ref{prop_1} and \ref{prop_2}, from the SUPG bracket in Definition \ref{def_SUPG_bracket}, we obtain a system of equations of the form
\begingroup
\addtolength{\jot}{2mm}
\begin{align}
&\left\langle \mathbf{w}, \mathbf{u}_t \right\rangle + \left\langle \mathbf{w}, \mathbf{q} \times \dd{H}{\mathbf{u}} \right\rangle - \left\langle \dd{H}{\rho}, \nabla \cdot \mathbf{w} \right\rangle \nonumber \\
&\hspace{2cm}- L\left(\frac{1}{\rho}\mathbf{w}, \theta; s\left(\tau, \mathbf{u}; T\right) + \tau S\left(\mathbf{u}; s\left(\tau, \mathbf{u}; T\right)\right)\right)  = 0 &\forall \mathbf{w} \in \mathbb{V}_u, \label{th_upw_u_eqn}\\
&\left\langle \phi, \rho_t \right\rangle + \left\langle \phi, \nabla \cdot \dd{H}{\mathbf{u}}\right\rangle = 0 &\forall \phi \in \mathbb{V}_\rho,\\
&\left\langle s + \tau S(\mathbf{u}; s), \theta_t \right\rangle + L\left(\frac{1}{\rho}\dd{H}{\mathbf{u}}, \theta; s + \tau S(\mathbf{u}; s)\right) = 0 &\forall s \in \mathbb{V}_\theta \label{th_upw_th_eqn}. 
\end{align}
\endgroup
Conversely, we find that the Poisson system can be recovered from these equations from the chain rule for $F$ as in \eqref{find_dFdt}. For any functional $F$, we have
\begin{align}
\frac{dF}{dt} &= \left\langle \dd{F}{\mathbf{u}}, \mathbf{u}_t \right\rangle + \left\langle \dd{F}{\rho}, \rho_t \right\rangle + \left\langle \dd{F}{\theta}, \theta_t \right\rangle \\
&= \left\langle \dd{F}{\mathbf{u}}, \mathbf{u}_t \right\rangle + \left\langle \dd{F}{\rho}, \rho_t \right\rangle + \left\langle s\left(\tau, \mathbf{u}; \dd{F}{\theta}\right) + S\left(\tau, \mathbf{u}; s\left(\tau, \mathbf{u}; \dd{F}{\theta}\right)\right), \theta_t \right\rangle = \{F, H\},
\end{align}
where the last equation follows directly from equations \eqref{th_upw_u_eqn} - \eqref{th_upw_th_eqn} and the previous proposition's observation regarding $T$ and $\dd{H}{\theta}$. Since the SUPG bracket is antisymmetric, we therefore find that equations \eqref{th_upw_u_eqn} - \eqref{th_upw_th_eqn} are still energy conserving.\\ \\
Finally, we recall that the SUPG method aims at dissipating $\|\theta\|_2^2$ in a manner proportional to the thermal field transport term. Since energy conservation was re-established by adding a corresponding antisymmetric, sign-opposite term to the bracket, this may raise concerns on whether the latter term has an anti-dissipative effect. However, the effect of the antisymmetric term is to transfer between potential and kinetic energy in order to conserve the total energy. Specifically, for the latter two sub-energies, given by
\begin{equation}
K = \frac{1}{2} \int_\Omega \rho |\mathbf{u}|^2 dx, \hspace{2cm} I =
\begin{cases}
c_v \rho \theta \pi \hspace{11mm} \text{(Euler equations)}, \\
\rho\theta\left(\frac{\rho}{2} + b\right) \;\;\; \text{(thermal SWE)},
\end{cases}
\end{equation}
the Poisson system \eqref{Poisson_system} yields
\begin{align}
F \!=\! K=&\;\; \rightarrow \;\; \frac{dK}{dt} = \dots + L\left(\mathbf{u}, \theta; s\left(\tau, \mathbf{u}; \dd{H}{\theta}\right) + S\left(\tau, \mathbf{u}; s\left(\tau, \mathbf{u}; \dd{H}{\theta}\right)\right)\right)\\
F = I=&\;\; \rightarrow \;\; \frac{dI}{dt} \;\;= \dots - L\left(\mathbf{u}, \theta; s\left(\tau, \mathbf{u}; \dd{I}{\theta}\right) \;+ S\left(\tau, \mathbf{u}; s\left(\tau, \mathbf{u}; \dd{I}{\theta}\;\right)\right)\right),
\end{align}
noting that $\dd{H}{\theta} = \dd{I}{\theta}$. The antisymmetric term is a rather complicated term linking the velocity and potential temperature field, and diffusive mixing of the latter field for the purpose of stabilisation can lead to potential energy increasing as well as decreasing. In particular, the antisymmetric term is not of the form of a positive definite stream-wise diffusion operator on the velocity field. \\ \\
Nonetheless, one might still have the concern that since the upwind stabilisation has the most significant effect at the grid-scale, the antisymmetric term would act as a source of grid-scale noise for the velocity. However, this is not necessarily the case, and for instance in \cite{natale2017scale} it was shown in numerical experiments for an energy conserving upwind discretisation of the 2D incompressible Euler equations, that the antisymmetric formulation causes energy to ``backscatter'' into low wavenumbers. This is possible because the extra antisymmetric components are nonlinear. To investigate whether the antisymmetric term is a source of high wavenumber noise, we consider in the numerical results section below semi-norms of the form
\begingroup
\addtolength{\jot}{4mm}
\begin{align}
& D\!G_\rho \coloneqq \sqrt{\sum_{K \in \mathcal{T}_h} \int_K \bigg( \Big(\pp{\rho}{x}\Big)^2 + \Big(\pp{\rho}{y}\Big)^2 \bigg) \; dx + \int_\Gamma \frac{1}{\Delta x_h} [\![\rho]\!]^2 \;dS}, \label{DG_measure_rho} \\
& D\!G_{\mathbf{u}} \coloneqq \sqrt{\int_\Omega \big((\nabla \cdot \mathbf{u})^2 + \omega^2 \big) \; dx}, \label{DG_measure_u}
\end{align}
\endgroup
for cells $K$ of a given tesselation $\mathcal{T}_h$, edge lengths $\Delta x_h$, and relative vorticity $\omega \in \mathbb{V}_q$, which is defined according to
\begin{align}
\langle \eta, \omega \rangle = - \left\langle \nabla^\perp \eta, \mathbf{u} \right\rangle + \left\langle\!\!\left\langle \eta, \mathbf{n}^\perp \cdot \mathbf{u} \right\rangle\!\!\right\rangle && \forall \eta \in \mathbb{V}_q.
\end{align}
The choice of $D\!G_\rho$ reflects semi-norms used in the stability analysis of DG methods (see e.g. \cite{arnold2002unified}). $D\!G_\mathbf{u}$ has been set up in a similar way, noting that for the div-conforming velocity field space, the jump term $[\![\mathbf{u}]\!]$ vanishes. In particular, we will find that the presence of the antisymmetric SUPG term does not lead to an increase of these semi-norms when compared to a discretisation that is not set up in an antisymmetric fashion (such as \eqref{nec_bracket} below). 
\begin{remark} \label{urho_upw}
To also incorporate energy conserving upwind stabilisation methods for the transport terms for $\mathbf{u}$ and $\rho$ in the thermodynamic equation set \eqref{th_upw_u_eqn} - \eqref{th_upw_th_eqn}, we can use the framework provided in \cite{wimmer2020energy}, where a bracket for the rotating shallow water equations was considered. The bracket is identical to the one introduced above, up to excluding the two thermal terms (i.e. setting $\dd{H}{\theta} = 0$ and $\dd{F}{\theta} = 0$ in $\{F, H\}$ in Definition \ref{def_SUPG_bracket}), and we can directly apply the transport methods here. The shallow water upwind stabilised bracket is given by
\begingroup
\addtolength{\jot}{1em}
\begin{align}
\{F, H\} \coloneqq&  \scalebox{1.0}[1]{$\left\langle \nabla^\perp \left(\rho \mathbb{U}\left(\rho, \dd{F}{\mathbf{u}}\right) \cdot \mathbb{U}\left(\rho, \dd{H}{\mathbf{u}}\right)^\perp \right), \mathbf{u} \right\rangle \!-\! \int_\Gamma \left[\!\left[\rho \mathbb{U}\left(\rho, \dd{F}{\mathbf{u}}\right) \cdot \mathbb{U}\left(\rho, \dd{H}{\mathbf{u}}\right)^\perp\right]\!\right] \mathbf{n}^\perp \cdot \tilde{\mathbf{u}} \; dS$} \label{var_Gv2_bracket_u_a}\\
& \scalebox{1.0}[1]{$-\! \left\langle \rho \; \mathbb{U}\left(\rho, \frac{\delta F}{\delta \mathbf{u}}\right), \nabla \frac{\delta H}{\delta \rho}\right\rangle \!+\! \int_\Gamma \left[\!\left[\frac{\delta H}{\delta \rho} \mathbb{U}\left(\rho, \frac{\delta F}{\delta \mathbf{u}}\right)\right]\!\right] \tilde{\rho} \; dS \!-\! \left\langle \rho \mathbb{U}\left(\rho, \dd{F}{\mathbf{u}}\right), f \mathbb{U}\left(\rho, \dd{H}{\mathbf{u}}\right)^\perp\right\rangle$} \label{var_Gv2_bracket_u_f}\\
& \scalebox{1.0}[1]{$+\!  \left\langle \rho \; \mathbb{U}\left(\rho, \frac{\delta H}{\delta \mathbf{u}}\right), \nabla \frac{\delta F}{\delta \rho}\right\rangle \!-\! \int_\Gamma \left[\!\left[\frac{\delta F}{\delta \rho} \mathbb{U}\left(\rho, \frac{\delta H}{\delta \mathbf{u}}\right)\right]\!\right] \tilde{\rho} \; dS$}, \label{var_Gv2_bracket_D_a}
\end{align}
\endgroup
where \eqref{var_Gv2_bracket_u_a} corresponds to the curl part of the velocity transport term's vector-invariant form (see \eqref{vector_invariant}), \eqref{var_Gv2_bracket_u_f} to the remaining terms in the momentum equation, and \eqref{var_Gv2_bracket_D_a} to density transport. As in \eqref{theta_SUPG_eqn_CP}, the gradient operators are considered element-wise as they are applied to functions that may be discontinuous across facets. For the purpose of density transport stabilisation, a velocity recovery operator $\mathbb{U}$ of the form
\begin{equation}
\hspace{2cm} \mathbb{U}(\rho, \mathbf{m}) \colon  \mathbb{V}_\rho \times \mathbb{V}_u \longrightarrow \mathbb{V}_u \; \; \; \; \text{such that} \; \; \; \; \left\langle \rho \mathbf{v}, \mathbb{U} \right\rangle = \left\langle \mathbf{v}, \mathbf{m} \right\rangle \hspace{15mm} \forall \mathbf{v} \in \mathbb{V}_u, \label{bbU_eqn}
\end{equation}
was introduced (for details, see \cite{wimmer2020energy}). Note that $\mathbb{U}$ corresponds to a discrete division by $\rho$, and recalling that $\dd{H}{\mathbf{u}}$ is given by the discrete flux $P_{\mathbb{V}_u}(\rho\mathbf{u})$, we find that $\mathbb{U}\left(\rho, \dd{H}{\mathbf{u}}\right) \in \mathbb{V}_u$ corresponds to the advecting velocity in the momentum and continuity equations. The bracket then leads to an upwind stabilised density equation of the form
\begin{align}
\left\langle \phi, \rho_t \right\rangle = \left\langle \rho \; \mathbb{U}\left(\rho, \frac{\delta H}{\delta \mathbf{u}}\right), \nabla \phi \right\rangle - \int_\Gamma \left[\!\left[\phi \mathbb{U}\left(\rho, \frac{\delta H}{\delta \mathbf{u}}\right)\right]\!\right] \tilde{\rho} \; dS && \forall \phi \in \mathbb{V}_\rho,
\end{align}
which corresponds to the standard DG upwind method. For the velocity equation, we obtain
\begin{align}
\left\langle \mathbf{w}, \mathbf{u}_t \right\rangle &= \left\langle \nabla^\perp \left(\rho \mathbb{U}(\rho, \mathbf{w}) \cdot \mathbb{U}\left(\rho, \dd{H}{\mathbf{u}}\right)^\perp \right), \mathbf{u} \right\rangle - \int_\Gamma \left[\!\left[\rho \mathbb{U}(\rho, \mathbf{w}) \cdot \mathbb{U}\left(\rho, \dd{H}{\mathbf{u}}\right)^\perp\right]\!\right] \mathbf{n}^\perp \cdot \tilde{\mathbf{u}} \; dS \nonumber \\
& - \left\langle \rho \mathbb{U}(\rho, \mathbf{w}), f \mathbb{U}\left(\rho, \dd{H}{\mathbf{u}}\right)^\perp\right\rangle \nonumber \\
&- \left\langle \rho \; \mathbb{U}(\rho, \mathbf{w}), \nabla \frac{\delta H}{\delta \rho}\right\rangle + \int_\Gamma \left[\!\left[\frac{\delta H}{\delta \rho} \mathbb{U}(\rho, \mathbf{w})\right]\!\right] \tilde{\rho} \; dS  \hspace{2cm} \forall \mathbf{w} \in \mathbb{V}_u,
\end{align}
where upwind stabilisation was applied to the curl part of the velocity transport term's vector-invariant form, i.e. $(\nabla^\perp \cdot \mathbf{u})\mathbf{u}^\perp$ (for details, see \cite{natale2016variational}). \\ \\
Note that the velocity recovery operator is applied to all test functions on the right hand side of the momentum equation. For a consistent use of test functions, when extending this bracket to the thermal case, we then also apply this operator to the test functions of the corresponding thermal bracket term (third term in \eqref{SUPG_bracket}). To maintain antisymmetry, we then also need to apply it in the thermal field transport term (fifth term in \eqref{SUPG_bracket}), and we arrive at a fully upwinded bracket given by \eqref{var_Gv2_bracket_u_a} - \eqref{var_Gv2_bracket_D_a}, together with
\begin{align}
\begin{split}
&+ L\left(\mathbb{U}\left(\rho, \dd{F}{\mathbf{u}}\right), \theta; s\left(\tau, \mathbf{u}; \dd{H}{\theta}\right) +\tau S\left(\mathbf{u}; s\left(\tau, \mathbf{u}; \dd{H}{\theta}\right)\right)\right) \\
&- L\left(\mathbb{U}\left(\rho, \dd{H}{\mathbf{u}}\right), \theta; s\left(\tau, \mathbf{u}; \dd{F}{\theta}\right) + \tau S\left(\mathbf{u}; s\left(\tau, \mathbf{u}; \dd{F}{\theta}\right)\right)\right). \label{bbU_SUPG_bracket}
\end{split}
\end{align}
The flux recovered velocity $\mathbb{U}\left(\rho, \dd{H}{\mathbf{u}}\right)$ then also serves as the advecting velocity for the thermal field transport equation, given by
\begin{align}
\left\langle s + \tau S(\mathbf{u}; s), \theta_t \right\rangle + L\left(\mathbb{U}\left(\rho, \dd{H}{\mathbf{u}}\right), \theta; s + \tau S(\mathbf{u}; s)\right) = 0 && \forall s \in \mathbb{V}_\theta. \label{bbU_SUPG_theta}
\end{align}
Finally, we note that the 2D velocity transport and Coriolis terms in \eqref{var_Gv2_bracket_u_a} and \eqref{var_Gv2_bracket_u_f} above can be extended readily to the three-dimensional case \cite{natale2016compatible}, which in our case leads to
\begin{align}
\begin{split}
 &\left\langle \nabla \times \left(\rho \mathbb{U}\left(\rho, \dd{F}{\mathbf{u}}\right) \times \mathbb{U}\left(\rho, \dd{H}{\mathbf{u}}\right) \right), \mathbf{u} \right\rangle - \int_\Gamma \left\{\!\left\{\mathbf{n} \times \left(\rho \mathbb{U}\left(\rho, \dd{F}{\mathbf{u}}\right) \times \mathbb{U}\left(\rho, \dd{H}{\mathbf{u}}\right)\right)\right\}\!\right\} \cdot \tilde{\mathbf{u}} \; dS \\
& - \left\langle \rho \mathbb{U}\left(\rho, \dd{F}{\mathbf{u}}\right), 2 \mathbf{\Omega} \times \mathbb{U}\left(\rho, \dd{H}{\mathbf{u}}\right)\right\rangle, \label{3D_u_advection_coriolis}
\end{split}
\end{align}
where
\begin{equation}
\{\!\{\mathbf{n} \times \mathbf{w} \}\!\} = \mathbf{n}^+ \times \mathbf{w}^+ + \mathbf{n}^- \times \mathbf{w}^-.
\end{equation}
\end{remark}
\section{Time discretisation} \label{section_time_discretisation}
In this section, we review the energy conserving time discretisation applied to \eqref{th_upw_u_eqn} - \eqref{th_upw_th_eqn} in order to obtain a fully energy conserving scheme, and describe the resulting non-linear set of equations. Further, we include a brief discussion on an approximately energy conserving time discretisation applied to the fully upwind stabilised setup as given in Remark \ref{urho_upw}.
\subsection{Poisson integrator}
To confirm the energy conserving property of the space discretised equations \eqref{th_upw_u_eqn} - \eqref{th_upw_th_eqn}, we apply an energy conserving time discretisation, thus expecting energy conservation up to solver tolerance. It is given by a Poisson integrator as introduced in \cite{cohen2011linear}, and can be applied to the framework used here as detailed in \cite{wimmer2020energy}. For the prognostic fields $\mathbf{z} = (\mathbf{u}, \rho, \theta)$, we have
\begin{equation}
\mathbf{z}^{n+1} = \mathbf{z}^n + \Delta t J\Big(\frac{\mathbf{z}^{n+1} + \mathbf{z}^n}{2}\Big) \Big(\overline{\dd{H}{\mathbf{u}}}, \overline{\dd{H}{\rho}}, \overline{\dd{H}{\theta}}\Big)^T,\label{J_system_discretised}
\end{equation}
where the skew symmetric transformation $J$ is related to the almost Poisson bracket via
\begin{equation}
\{F, H\} = \left\langle \dd{F}{\mathbf{z}},  J(\mathbf{z}) \dd{H}{\mathbf{z}} \right\rangle,
\end{equation}
and the time averaged Hamiltonians are given by
\begin{equation}
\overline{\dd{H}{\mathbf{u}}} \coloneqq \int_0^1 \dd{}{\mathbf{u}} H(\mathbf{z}^n + s(\mathbf{z}^{n+1} - \mathbf{z}^n)) ds, \label{dH_bar}
\end{equation}
and similarly for the variations in $\rho$ and $\theta$. The expressions can be integrated exactly for the thermal shallow water Hamiltonian, leading to
\begin{equation}
\overline{\dd{H}{\mathbf{u}}} = \frac{1}{3} P_{\mathbb{V}_u} \big(\rho^n \mathbf{u}^n + \frac{1}{2} \rho^n \mathbf{u}^{n+1} + \frac{1}{2} \rho^{n+1} \mathbf{u}^n + \rho^{n+1} \mathbf{u}^{n+1} \big), \label{dHdu_bar}
\end{equation}
and similarly for the other two variations. However, for the Euler Hamiltonian, we find that the internal energy is a non-polynomial function in $\rho$ and $\theta$, thus requiring an approximate integration of \eqref{dH_bar} for the variations in $\rho$ and $\theta$. For the numerical results below, a fourth order Gaussian quadrature was used.\\ \\
The resulting nonlinear system of equations is given by
\begingroup
\addtolength{\jot}{3mm}
\begin{align}
&\left\langle \mathbf{w}, \mathbf{u}^{n+1} - \mathbf{u}^n \right\rangle + \Delta t\Bigg(\left\langle \mathbf{w}, \bar{\mathbf{q}} \times \overline{\dd{H}{\mathbf{u}}} \right\rangle - \left\langle \overline{\dd{H}{\rho}}, \nabla \cdot \mathbf{w} \right\rangle \nonumber \\
&\hspace{4cm}- L\left(\frac{1}{\bar{\rho}}\mathbf{w}, \theta; s\left(\tau, \bar{\mathbf{u}}; \bar{T}\right) + \tau S\left(\bar{\mathbf{u}}; s\left(\tau, \bar{\mathbf{u}}; \bar{T}\right)\right)\right) \Bigg) = 0 & \forall \mathbf{w} \in \mathbb{V}_u, \label{th_upw_u_eqn_dscr} \\
&\left\langle \phi, \rho^{n+1} - \rho^n \right\rangle + \Delta t\left\langle \phi, \nabla \cdot \overline{\dd{H}{\mathbf{u}}}\right\rangle = 0 & \forall \phi \in \mathbb{V}_\rho,\\
&\langle s + \tau S(\bar{\mathbf{u}}; s), \theta^{n+1} - \theta^n \rangle + \Delta t L\left(\frac{1}{\bar{\rho}}\overline{\dd{H}{\mathbf{u}}}, \bar{\theta}; s + \tau S(\bar{\mathbf{u}}; s)\right) = 0 & \forall s \in \mathbb{V}_\theta, \label{th_upw_th_eqn_dscr}
\end{align}
\endgroup
for midpoint time average $\bar{\rho} = (\rho^n + \rho^{n+1})/2$, and similarly for $\bar{\mathbf{u}}$ and $\bar{\mathbf{q}}$, with $\mathbf{q}^n$ and $\mathbf{q}^{n+1}$ given by the diagnostic scalar vorticity equation \eqref{discr_diag_q}, using the prognostic fields at time levels $n$ and $n+1$, respectively. Further, $\bar{T}$ is defined analogously to $T$, i.e. such that $\overline{\dd{H}{\theta}} = P_{\mathbb{V}_\theta}(\bar{T})$. Additionally, we remark that the derivation of the fully discretised thermal field equation \eqref{th_upw_th_eqn_dscr} follows the proof of Proposition \ref{prop_1} in a time discretised form, with $\theta_t$ replaced by $(\theta^{n+1} - \theta^n)/\Delta t \in \mathbb{V}_\theta$. \\ \\
To solve for $\mathbf{z}^{n+1}$, given the fully discretised residual $\mathbf{R}(\mathbf{z}^{n+1})$ above (left-hand sides of \eqref{th_upw_u_eqn_dscr} - \eqref{th_upw_th_eqn_dscr}), we apply the same procedure as detailed in \cite{wimmer2020energy} and revert to a Picard iteration method. For update $\delta \mathbf{z} = \mathbf{z}^{n+1, k+1} - \mathbf{z}^{n+1, k}$ with unknown next time step estimate $\mathbf{z}^{n+1, k+1}$, we set
\begin{equation}
\dd{\mathbf{R}'}{\mathbf{z}}(\delta \mathbf{z}) = - \mathbf{R}(\mathbf{z}^{n+1, k}), \label{Picard_iteration_method}
\end{equation}
and $\mathbf{z}^{n+1,0} = \mathbf{z}^n$. The left hand side corresponds to the Jacobian of a linearised version of $\mathbf{R}$ without Hamiltonian projections, and for the Euler equations it is given by
\begin{equation}
\hspace{-2mm}\dd{\mathbf{R}'}{\mathbf{z}}(\delta \mathbf{z}) = \begin{pmatrix*}[l]
\left\langle \delta \mathbf{u}, \mathbf{w} \right\rangle \!+\! \frac{\Delta t}{2} \Big(\left\langle2 \mathbf{\Omega} \!\times\! \delta\mathbf{u}, \mathbf{w} \right\rangle \! - \! \left\langle g \delta \rho + c_p (\bar{\bar{\theta}}\delta \pi + \delta \theta \bar{\bar{\pi}}), \nabla \cdot \mathbf{w} \right\rangle \!-\! c_p \left\langle\bar{\bar{\pi}} \nabla \delta \theta, \mathbf{w} \right\rangle\Big)\\
\left\langle \delta \rho, \phi \right\rangle + \frac{\Delta t}{2} \left\langle \bar{\bar{\rho}} \nabla \cdot \mathbf{\delta u}, \phi \right\rangle \\
\left\langle \delta \theta, \gamma \right\rangle  \label{Linearised_Jacobian}
\end{pmatrix*}\!,
\end{equation}
where double-barred entries correspond to background fields, with $\bar{\bar{\pi}} = \pi(\bar{\bar{\rho}}, \bar{\bar{\theta}})$, $\bar{\bar{\mathbf{u}}} = \mathbf{0}$, and $\delta \pi = \pp{\pi}{\rho}(\bar{\bar{\rho}}, \bar{\bar{\theta}})\delta \rho + \pp{\pi}{\theta}(\bar{\bar{\rho}}, \bar{\bar{\theta}})\delta \theta$. Note that in view of the test cases to follow, we assumed the background thermal field to be constant (equal to $g$ for the thermal shallow water equations, and an isentropic background potential temperature for the Euler equations). This leads to a vanishing $\nabla \bar{\bar{\theta}}$ term in the velocity and thermal field equation, respectively, thus uncoupling the latter from the density and velocity equations. Given the right-hand side $- \mathbf{R}(\mathbf{z}^{n+1, k})$, we can then first solve for $\delta \theta$, followed by a mixed solve for the density and velocity updates. Note that since this paper focuses on the underlying space discretisation, the nonlinear solve procedure is kept simple. The discretisation can in principle be applied equally to fully three-dimensional problems and scenarios with a non-constant potential background temperature (by approximately eliminating $\delta \theta$ in the linearised momentum equation). A suitable nonlinear solver strategy for this case can be found in \cite{bendall2019compatible}.
\subsection{Approximately energy conserving scheme} \label{section_approx_ec}
While both the space and time discretisations discussed above conserve energy, the underlying Picard iteration method does not. In the numerical tests below, we find that the number of Picard iterations required to achieve energy conservation up to machine precision is much higher than is usually considered in forecasting models. This suggests that we may simplify the fully discretised scheme, ensuring that the additional energy error is smaller than the error due to a small number of Picard iterations. In particular, such a simplification was considered in \cite{wimmer2020energy} (Remark 4) for the fully upwind stabilised equations resulting from the bracket \eqref{var_Gv2_bracket_u_a} - \eqref{var_Gv2_bracket_D_a}. In order to avoid a computation of $\overline{\mathbb{U}} \coloneqq \mathbb{U}\left(\bar{\rho}, \overline{\dd{H}{\mathbf{u}}}\right)$, i.e.
\begin{align}
\left\langle \frac{1}{2}(\rho^n + \rho^{n+1}) \mathbf{v}, \bar{\mathbb{U}} \right\rangle = \left\langle \mathbf{v},  \frac{1}{3}\big(\rho^n \mathbf{u}^n + \frac{1}{2} \rho^n \mathbf{u}^{n+1} + \frac{1}{2} \rho^{n+1} \mathbf{u}^n + \rho^{n+1} \mathbf{u}^{n+1} \big) \right\rangle && \forall \mathbf{v} \in \mathbb{V}_u,
\end{align}
the Hamiltonian variation in $\mathbf{u}$ is instead time-discretised according to
\begin{equation}
\dd{H}{\mathbf{u}} = P_{\mathbb{V}_u}(\rho \mathbf{u}) \; \; \; \rightarrow \; \; \; P_{\mathbb{V}_u}(\bar{\rho}\bar{\mathbf{u}}), \label{midpoint_dHdu}
\end{equation}
which leads to
\begin{equation}
\overline{\mathbb{U}} = \bar{\mathbf{u}},
\end{equation}
point-wise. The deviation \eqref{midpoint_dHdu} to the Poisson integrator (which requires \eqref{dHdu_bar}) can equally be applied when the upwind stabilised bracket \eqref{var_Gv2_bracket_u_a} - \eqref{var_Gv2_bracket_D_a} is extended to include a thermal part of the form \eqref{bbU_SUPG_bracket}. The resulting fully discretised scheme, which is energy conserving up to the time-discretisation of $\dd{H}{\mathbf{u}}$, is then given by
\begingroup
\allowdisplaybreaks
\addtolength{\jot}{4mm}
\begin{align}
&\left\langle \bar{\rho} \mathbf{w}, \mathbf{u}^{n+1} - \mathbf{u}^n \right\rangle = \Delta t \bigg( \left\langle \nabla^\perp \left(\bar{\rho} \mathbf{w} \cdot \bar{\mathbf{u}}^\perp \right), \bar{\mathbf{u}} \right\rangle - \int_\Gamma \left[\!\left[\bar{\rho} \mathbf{w} \cdot \bar{\mathbf{u}}^\perp\right]\!\right] \mathbf{n}^\perp \cdot \tilde{\bar{\mathbf{u}}} \; dS \nonumber \\
& \hspace{42mm} - \left\langle \bar{\rho} \mathbf{w}, f \bar{\mathbf{u}}^\perp\right\rangle - \left\langle \bar{\rho} \mathbf{w}, \nabla \overline{\frac{\delta H}{\delta \rho}}\right\rangle + \int_\Gamma \left[\!\left[\overline{\frac{\delta H}{\delta \rho}} \mathbf{w}\right]\!\right] \tilde{\bar{\rho}} \; dS \nonumber \\
& \hspace{42mm} + L\left(\bar{\rho}\mathbf{w}, \bar{\theta}; s(\tau, \bar{\mathbf{u}}, \bar{T}) + \tau S(\bar{\mathbf{u}}; s(\tau, \bar{\mathbf{u}}, \bar{T}))\right) \bigg)  & \forall \mathbf{w} \in \mathbb{V}_u, \label{u_eqn_mid} \\
&\left\langle \phi, \rho^{n+1} - \rho^{n} \right\rangle = \Delta t \left(\left\langle \bar{\rho} \bar{\mathbf{u}}, \nabla \phi \right\rangle - \int_\Gamma \left[\!\left[\phi \bar{\mathbf{u}}\right]\!\right] \tilde{\bar{\rho}} \; dS \right) & \forall \phi \in \mathbb{V}_\rho, \label{rho_eqn_mid} \\
&\left\langle s + \tau S(\bar{\mathbf{u}}; s), \theta^{n+1} - \theta^n \right\rangle = - \Delta t L\left(\bar{\mathbf{u}}, \bar{\theta}; s + \tau S(\bar{\mathbf{u}}; s)\right) & \forall s \in \mathbb{V}_\theta, \label{theta_eqn_mid}
\end{align}
\endgroup
with $L$, $S$, $s$, and $T$ as described in Definitions \ref{def_gamma_upw_L_operator} and \ref{definition_s} as well as Proposition \ref{prop_2}, respectively. Note that \eqref{rho_eqn_mid} and \eqref{theta_eqn_mid} correspond to the DG upwind and SUPG methods, respectively, together with an implicit midpoint time discretisation. Further, the first two terms on the right-hand side of the momentum equation \eqref{u_eqn_mid} correspond to an upwind stabilisation as introduced in \cite{natale2016variational}, again together with an implicit midpoint time discretisation. Additionally, the last two lines of the momentum equation correspond to the Coriolis term as well as other forcing terms that are formulated according to the almost Poisson bracket's underlying antisymmetry. Finally, the momentum equation's test functions are weighted by $\bar{\rho}$ to avoid the occurrence of $\mathbb{U}(\bar{\rho}, \mathbf{w})$, i.e. the velocity recovery operator applied to test functions (for details, see \cite{wimmer2020energy}).\\ \\
Altogether, the resulting scheme resembles discretisations derived outside the Hamiltonian framework, such as the second one introduced in Section \ref{section_nec} below, up to the weighted test functions and modified forcing terms in the momentum equation. In particular, this facilitates the implementation of this scheme, since forcing terms in the momentum equation are often handled separately in code bases. Finally, from the point of view of computational cost, we find that in the approximately energy conserving scheme \eqref{u_eqn_mid} - \eqref{theta_eqn_mid}, we are additionally required to compute $\overline{\dd{H}{\rho}} \in \mathbb{V}_\rho$ and $s(\tau, \bar{\mathbf{u}}, T) \in \mathbb{V}_\theta$ in each Picard iteration. For the Euler equations, this corresponds to a projection into a DG space, together with a computation in a space that is discontinuous in the horizontal direction, which overall leads to a only small increase in computational cost when compared to the second scheme in Section \ref{section_nec}.
\section{Numerical results} \label{section_numerical_results}
In this section, we confirm numerically the upwind-stabilised and energy conserving properties of the fully discretised schemes \eqref{th_upw_u_eqn_dscr} - \eqref{th_upw_th_eqn_dscr} and \eqref{u_eqn_mid} - \eqref{theta_eqn_mid}. For comparison purposes, we first describe two non-energy conserving space discretisations. We then consider test cases consisting of a (perturbed) steady state scenario for the thermal shallow water equations, as well as a cold and hot air bubble scenario for the Euler equations.
\subsection{Comparison to non-energy conserving space discretisations} \label{section_nec}
The first non-energy conserving space discretisation is formulated such that the Poisson time integrator can be applied to it, allowing for a comparison that focuses on the exact energy conserving properties of the SUPG stabilised space discretisation. The second is formulated using a standard treatment of the Euler equations, and is used in comparison with the fully upwind stabilised scheme as described in Section \ref{section_approx_ec}, with a smaller number of Picard iterations in the simulation runs.\\ \\
For the first discretisation, we use a non-antisymmetric bracket of the form
\begingroup
\allowdisplaybreaks
\addtolength{\jot}{4mm}
\begin{align}
\begin{split}
\hspace{-3mm} \{F, H\} = &- \left\langle \dd{F}{\mathbf{u}}, \mathbf{q} \times \dd{H}{\mathbf{u}} \right\rangle \\
& + \left\langle \dd{H}{\rho}, \nabla \cdot \dd{F}{\mathbf{u}} \right\rangle - \left\langle \nabla \cdot \Big(\frac{1}{\rho} \dd{H}{\theta} \dd{F}{\mathbf{u}}\Big), \theta \right\rangle + \int_\Gamma \left[\!\left[\frac{1}{\rho} \dd{H}{\theta} \dd{F}{\mathbf{u}}\right]\!\right] \{\theta\} \; dS\\
& - \left\langle \dd{F}{\rho}, \nabla \cdot \dd{H}{\mathbf{u}} \right\rangle - L\left(\frac{1}{\rho}\dd{H}{\mathbf{u}}, \theta; s\left(\tau, \mathbf{u}; \dd{F}{\theta}\right) + S\left(\tau, \mathbf{u}; s\left(\tau, \mathbf{u}; \dd{F}{\theta}\right)\right)\right), \label{nec_bracket}
\end{split}
\end{align}
\endgroup\\
where $\{\cdot\}$ denotes the average across facets. Note that this bracket is antisymmetric with respect to the vorticity term and the density transport terms, but not with respect to the thermal field transport term. The last term is set to arrive at the SUPG stabilised form \eqref{SUPG_theta_ec} of the thermal field transport equation, while the third and fourth terms correspond to a simple stabilisation of the non-upwind stabilised thermal field bracket term:
\begin{equation}
\left\langle \frac{1}{\rho} \dd{H}{\theta} \nabla \theta, \dd{F}{\mathbf{u}} \right\rangle \; \; \; \rightarrow \; \; \; - \left\langle \nabla \cdot \Big(\frac{1}{\rho} \dd{H}{\theta} \dd{F}{\mathbf{u}}\Big), \theta \right\rangle + \int_\Gamma \left[\!\left[\frac{1}{\rho} \dd{H}{\theta} \dd{F}{\mathbf{u}}\right]\!\right] \{\theta\} \; dS.
\end{equation}
This allows us to investigate the loss in energy due to a non-energy conserving implementation of the SUPG method, noting that we couple the non-antisymmetric bracket with the fully energy conserving Poisson integrator.\\ \\
The second discretisation is derived directly from the continuous Euler equations \eqref{Euler_u_eqn} - \eqref{Euler_th_eqn}, disregarding the Hamiltonian framework. For this purpose, we use the same types of upwind stabilisation as discussed in Section \ref{section_formulation}, i.e. standard DG upwinding for $\rho$, SUPG for $\theta$, and upwinding for the curl part of the velocity transport term's vector-invariant form. The resulting weak form is then given by
\begingroup
\allowdisplaybreaks
\addtolength{\jot}{4mm}
\begin{align}
&\left\langle \mathbf{w}, \mathbf{u}^{n+1} - \mathbf{u}^n \right\rangle = \Delta t \bigg( \left\langle \nabla^\perp \left(\mathbf{w} \cdot \bar{\mathbf{u}}^\perp \right), \bar{\mathbf{u}} \right\rangle - \int_\Gamma \left[\!\left[\mathbf{w} \cdot \bar{\mathbf{u}}^\perp\right]\!\right] \mathbf{n}^\perp \cdot \tilde{\bar{\mathbf{u}}} \; dS \nonumber \\
&\hspace{40mm} - \left\langle \nabla \cdot \mathbf{w}, |\bar{\mathbf{u}}|^2 \right\rangle + \left\langle \mathbf{w}, g\mathbf{k} \right\rangle\nonumber \\
& \hspace{40mm}- c_p \left\langle\nabla \cdot (\bar{\theta} \mathbf{w}), \bar{\pi} \right\rangle + c_p \int_\Gamma \left[\!\left[\bar{\theta} \mathbf{w}\right]\!\right] \{\bar{\pi}\} \; dS \bigg) = 0 &\forall \mathbf{w} \in \mathbb{V}_u,  \label{u_eqn_nec_mid}\\
&\left\langle \phi, \rho^{n+1} - \rho^{n} \right\rangle = \Delta t \left(\left\langle \bar{\rho} \bar{\mathbf{u}}, \nabla \phi \right\rangle - \int_\Gamma \left[\!\left[\phi \bar{\mathbf{u}}\right]\!\right] \tilde{\bar{\rho}} \; dS \right) & \forall \phi \in \mathbb{V}_\rho, \\
&\left\langle s + \tau S(\bar{\mathbf{u}}; s), \theta^{n+1} - \theta^n \right\rangle = - \Delta t L\left(\bar{\mathbf{u}}, \bar{\theta}; s + \tau S(\bar{\mathbf{u}}; s)\right) & \forall s \in \mathbb{V}_\theta, \label{theta_eqn_nec_mid}
\end{align}
\endgroup
for $\bar{\pi} = \pi(\bar{\rho}, \bar{\theta})$, and noting that for the Euler vertical slice test cases considered below, the Coriolis parameter $f$ is set zero. We use midpoint averages for the time discretisation, so that the resulting fully discretised density and thermal field equations are equal to the ones following from the approximately energy conserving discretisation described in Section \ref{section_approx_ec} above. In particular, in comparison to this non-Hamiltonian discretisation, we find that the adjustments used to derive the approximately energy conserving discretisation lie entirely in the momentum equation.
\subsection{Test cases}
Having formulated the Poisson bracket (Definition \ref{def_SUPG_bracket}) and the Poisson integrator, as well as the approximately energy conserving discretisation (Section \ref{section_approx_ec}) and two non-energy conserving schemes, we test for their energy-conserving properties and the qualitative field development. We first consider a thermal shallow water test case and an Euler test case, demonstrating energy conservation up to solver tolerance as well as an improved thermal field development when the SUPG method is applied. Given the shortcomings of the field development in the Euler test case, we then move on to the fully upwind stabilised setup, which we investigate using two Euler test cases. \\ \\
The mesh, finite element discretisation, and solver were implemented using the automated finite element toolkit Firedrake\footnote{for further details, see \cite{gmd-9-3803-2016, kirby2018solver, mcrae2016automated} or \url{http://firedrakeproject.org}} \cite{rathgeber2016firedrake}. The resulting systems of equations are solved using PETSc \cite{balay2019petsc, balay1997efficient}, noting that a hybridised solver (\cite{gibson2020slate} and e.g. \cite{shipton2018higher}) is used to solve for the mixed system in $(\delta \mathbf{u}, \; \delta \rho)$. Finally, for simplicity, we set the SUPG parameter to $\tau = \Delta t/2$ in all test cases; a further discussion on possible choices of $\tau$ can be found in \cite{kuzmin2010guide}.
\subsubsection{Energy conservation and thermal field upwinding}
\begin{figure}[ht]
\begin{center}
\includegraphics[width=1\textwidth]{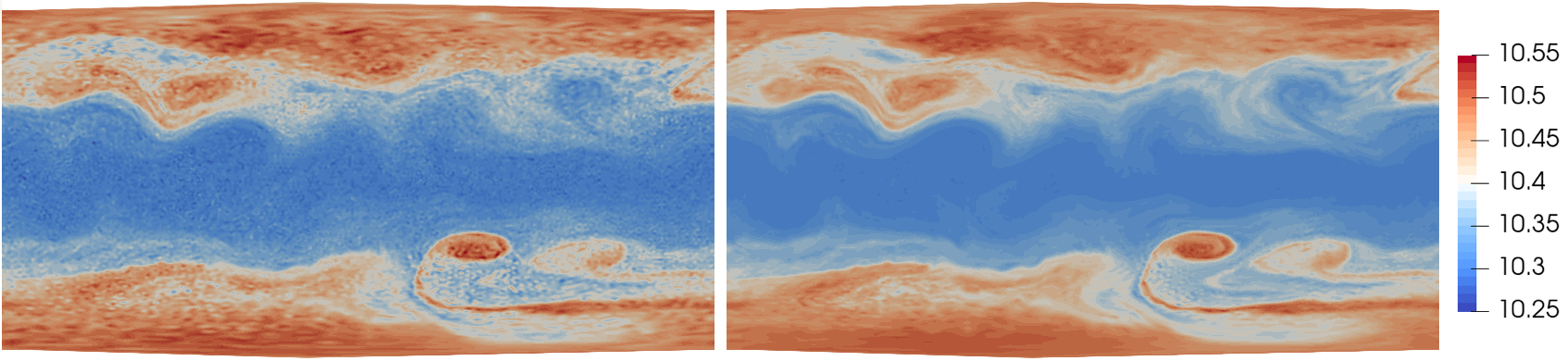}
\caption{Buoyancy fields for perturbed steady state test case after 50 days, for energy conserving discretisations. Left: no SUPG for buoyancy field; right: SUPG for buoyancy field. Mesh refinement level 4, $\Delta t = 8$ minutes, with 8 Picard iterations per time step.} \label{Thermal_W5_buoyancies}
\end{center}
\end{figure}

\begin{figure}[ht]
\begin{center}
\includegraphics[width=1\textwidth]{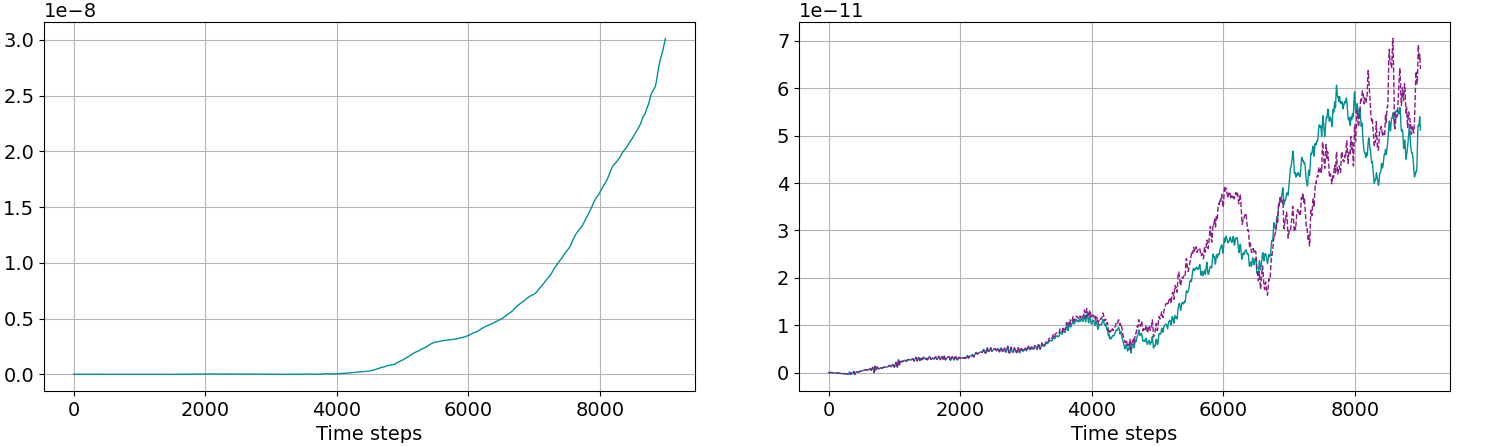}
\caption{Relative energy error development for perturbed steady state test case. Left: non-energy conserving bracket. Right: energy conserving bracket with SUPG for buoyancy (cyan) and energy conserving bracket without SUPG for buoyancy (dashed purple).} \label{Thermal_W5_energies}
\end{center}
\end{figure}

\begin{figure}[ht]
\begin{center}
\includegraphics[width=1\textwidth]{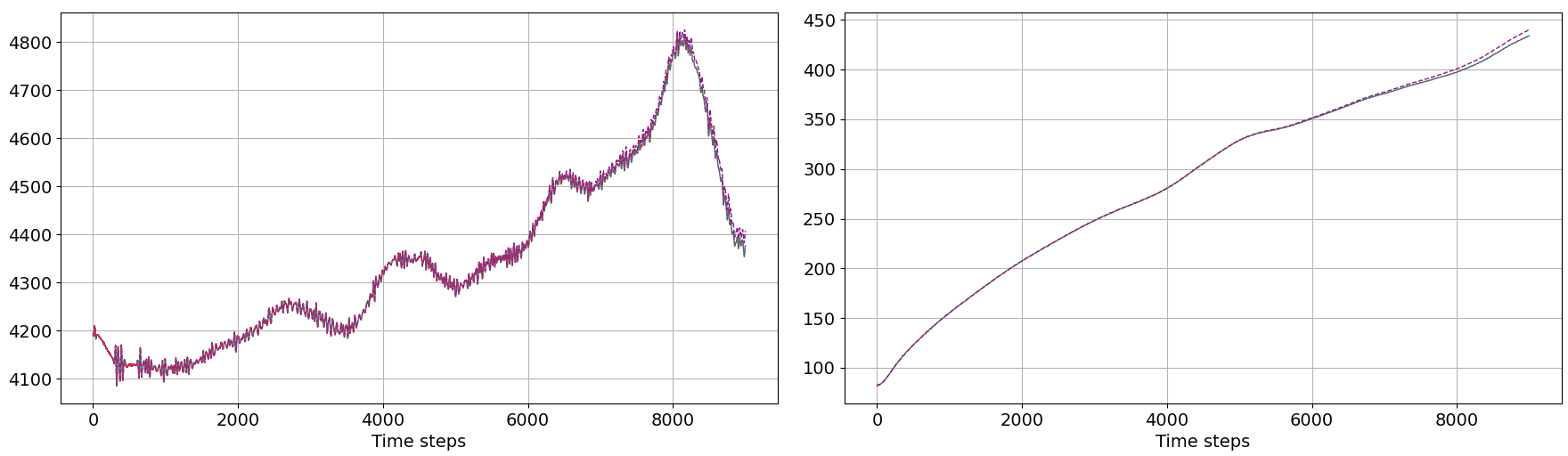}
\caption{Time development of $D\!G_\rho$ (left) and $D\!G_{\mathbf{u}}$ (right) for perturbed steady state test case. Cyan: energy conserving bracket with SUPG. Dashed purple: energy conserving bracket without SUPG. Dotted red: non-energy conserving bracket.} \label{W5_DG_rho_u}
\end{center}
\end{figure}

\begin{figure}[ht]
\begin{center}
\includegraphics[width=1\textwidth]{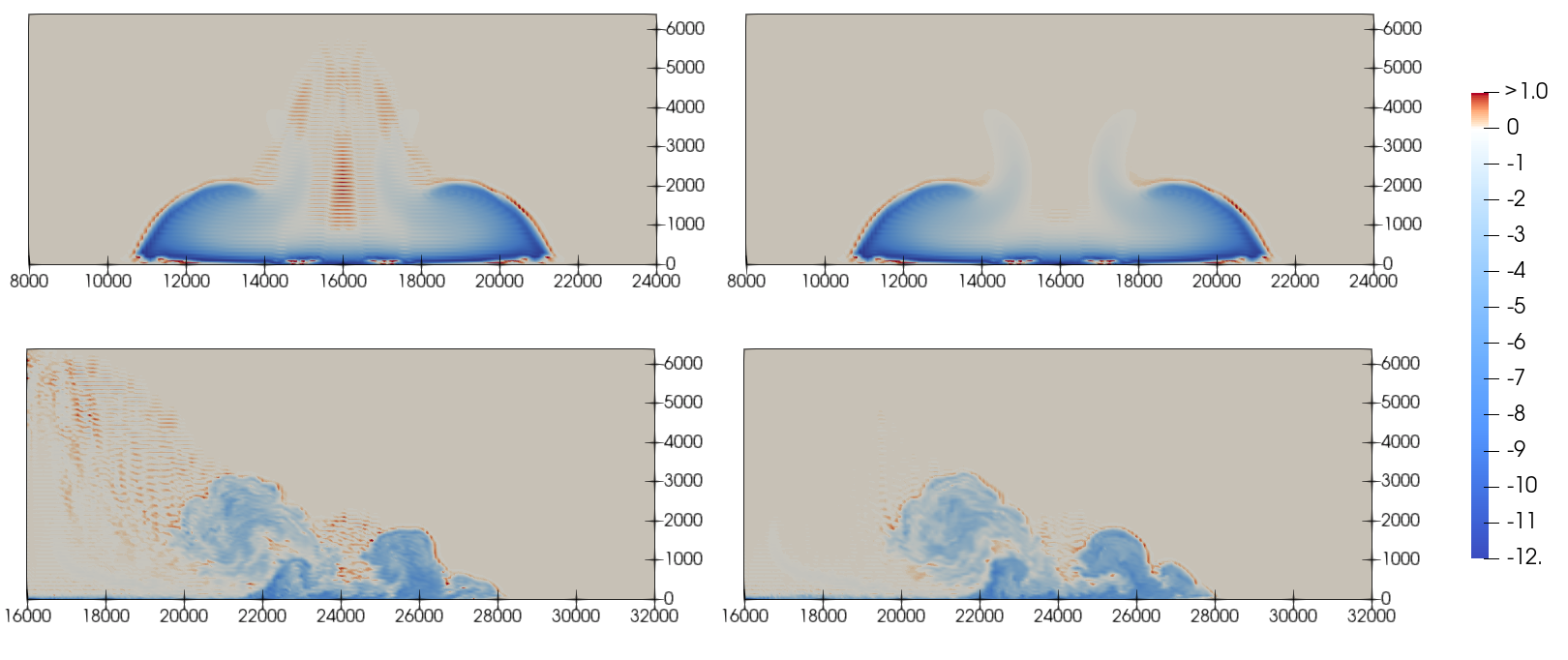}
\caption{Potential temperature fields for falling bubble test case after 400 s (top row) and 800 s (bottom row), for energy conserving discretisations. Left column: no SUPG for potential temperature field; right column: SUPG for potential temperature field. $\Delta x = \Delta z = 100$ m, $\Delta t = 0.5$ s, with 32 Picard iterations per time step.} \label{DB_potential_temperatures}
\end{center}
\end{figure}

\begin{figure}[ht]
\begin{center}
\includegraphics[width=1\textwidth]{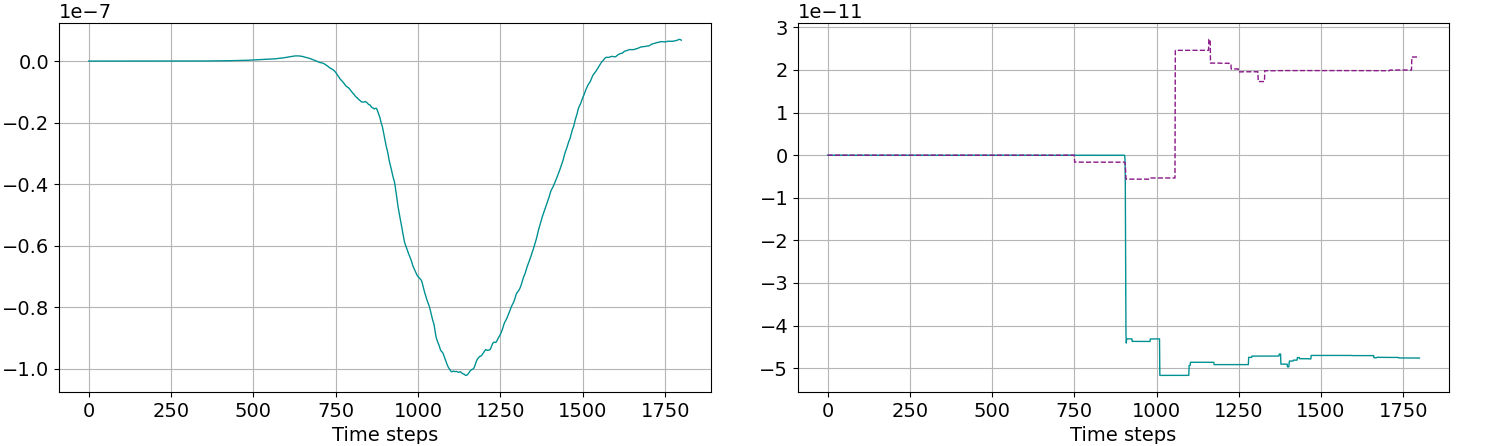}
\caption{Relative energy error development for falling bubble test case. Left: non-energy conserving bracket. Right: energy conserving bracket with SUPG for potential temperature (cyan) and energy conserving bracket without SUPG for potential temperature (dashed purple).} \label{DB_energies}
\end{center}
\end{figure}

\begin{figure}[ht]
\begin{center}
\includegraphics[width=1\textwidth]{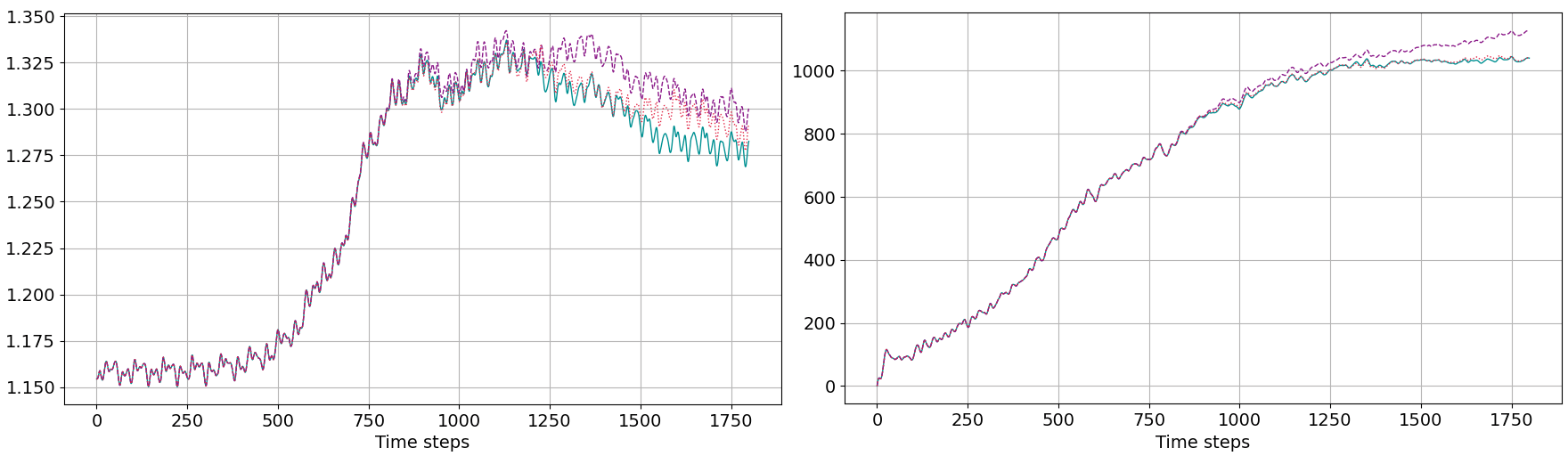}
\caption{Time development of $D\!G_\rho$ (left) and $D\!G_{\mathbf{u}}$ (right) for falling bubble test case. Cyan: energy conserving bracket with SUPG. Dashed purple: energy conserving bracket without SUPG. Dotted red: non-energy conserving bracket.} \label{DB_DG_rho_u}
\end{center}
\end{figure}

The thermal shallow water test case considered here is based on the fifth test case in \cite{williamson1992standard}, corresponding to a steady state spherical flow given by a geostrophic balance, which is perturbed by a mountain. Following \cite{ELDRED20191}, we set the thermal field such that the overall state is still balanced. The initial conditions are then given by
\begin{align}
&\mathbf{u} = u_0(-y, x, 0)/a,\\
&\rho = h - (a\Omega u_0+ u_0^2/2)\frac{z^2}{ga^2} - b,\\
&\theta =  g\bigg(1 + \epsilon \Big(\frac{\bar{\bar{\rho}}}{\rho}\Big)^2\bigg),
\end{align}
where $b= b_0(1 - r/R)$ describes the mountain's surface, for $R=\pi/9$, mountain height $b_0 = 2000$ m and $r = \text{min}(R, \sqrt{(\lambda - \lambda_c)^2 + (\theta - \theta_c)^2})$. $\lambda \in [-\pi, \pi]$ and $\theta \in [-\pi/2, \pi/2]$ denote longitude and latitude respectively, and the mountain's centre is chosen as $\lambda_c = -\pi/2$ and $\theta_c = \pi/6$. The sphere's radius, rotation rate, and the gravitational acceleration are given by $a=6371220$ m, $\Omega = 7.292 \times 10^{-5}$ s$^{-1}$, and $g = 9.810616$ ms$^{-2}$, respectively. Further, the mean height, wind speed, and thermal field parameter are set to $\bar{\bar{\rho}} = 5960$ m, $u_0 = 20$ m/s, and $\epsilon = 0.05$, respectively.\\ \\
The mesh is given by an icosahedral triangulation of the sphere and is built using a second (polynomial) order coordinate space. Refinement level 0 corresponds to 20 triangles, and for every higher level, each triangle is refined to 4 triangles (so that each increase corresponds to halving the cell side length $\Delta x$). In this test case the refinement level was set to $4$. The simulation is run for 50 days, with a time step of $\Delta t = 8$ minutes, and 8 Picard iterations for each time step. To also test for the convergence property of the energy conserving scheme including SUPG stabilisation, the simulation is additionally run as a steady state test case (i.e. without the perturbing mountain profile $b$), with a runtime of 50 days, time step of $\Delta t = 30$ minutes, 4 Picard iterations for each time step, and refinement levels $3$, $4$, and $5$. Note that this corresponds to the second test case in \cite{williamson1992standard} extended to include a thermal field.\\ \\
For the steady state, we find error ratios given by
\begin{center}
\begin{tabular}{|c|c|c|} 
 \hline
 $i$ & $e_\theta^i/e_\theta^{i+1}$ & $e_\mathbf{u}^i/e_\mathbf{u}^{i+1}$\\ 
 \hline
 3 & 4.00165 & 3.82207 \\
 4 & 4.00318 & 3.95878 \\
 \hline
\end{tabular}
\end{center}
\vspace{2mm}
where $e^i_\theta$ corresponds to the $L^2$ difference between the initial and final buoyancy fields at the $i^{th}$ refinement level, and similarly for $e^i_\mathbf{u}$. As expected from approximation theory for second order BDM elements, the convergence rate is of the order of $(\Delta x)^2$. For the perturbed setup, the resulting buoyancy fields for the discretised equations corresponding to the non-upwinded bracket \eqref{cts_Euler_bracket} and the SUPG bracket (Definition \ref{def_SUPG_bracket}) are depicted in Figure \ref{Thermal_W5_buoyancies}. Further, the relative energy error (i.e. $(E_t - E_0)/E_0$) development for the aforementioned brackets as well as the non-energy conserving one (i.e. \eqref{nec_bracket}) are depicted in Figure \ref{Thermal_W5_energies}. Finally, the time developments of $D\!G_\rho$ and $D\!G_{\mathbf{u}}$ are given in Figure \ref{W5_DG_rho_u}.\\
%
\\
\\
\\
As expected, we find that the incorporation of SUPG markedly reduces the occurrence of numerical oscillations in the buoyancy field development. Further, both energy conserving brackets lead to energy convergence up to solver tolerance, with an improvement by 3 orders of magnitude when compared to the non-energy conserving bracket. Additionally, while the error curves for the energy conserving schemes do not show a clear tendency across the simulation, the one for the non-energy conserving scheme does have one. While this test case was run for 50 simulated days only, climate simulations are typically run for many simulated years. In particular, while an energy error of the order of $10^{-8}$ may be seen as relatively small, the existing bias may aggravate the error substantially for such long scale simulations. Finally, we find that the development of $DG_\rho$ and $DG_\mathbf{u}$ is nearly identical for all three schemes, except for slightly larger values of $DG_\mathbf{u}$ for the energy-conserving setup without SUPG. The latter is likely due to noise in the thermal field affecting the velocity field. Altogether, this indicates that the antisymmetric forcing term corresponding to SUPG (last term in \eqref{th_upw_u_eqn}) does not lead to an increase in oscillations in the velocity and density fields.\\ \\
Next, we consider the brackets for the Euler equations. For this purpose, we use the next to lowest order density finite element space (k=2) and the corresponding spaces for the other fields as described in Section \ref{section_compatible_FEM}. Since the non-upwinded bracket \eqref{cts_Euler_bracket} leads to an unstable field development in the test case described below, we instead compare the SUPG bracket to one that also includes DG-upwinding for potential temperature in the horizontal direction, but no SUPG in the vertical. Note that such a bracket can be derived from the SUPG bracket by setting the SUPG parameter $\tau$ to zero. Further, note that the diagnostic vorticity equation \eqref{discr_diag_q} does not conserve vorticity near the boundary even in the absence of a baroclinic term, which may lead to an unbounded growth of enstrophy due to sources near the boundary (for more details, see \cite{BAUER2018171}). However, for the test case below this effect is small, and given our focus on thermal field transport methods we ignore it.\\ \\
The test case is given by a falling bubble in a vertical slice of the atmosphere, based on one as described in  \cite{straka1993numerical}. For this purpose, we consider a horizontally periodic rectangular domain $\Omega$ of $32$ km length and $6.4$ km height, with a constant potential temperature background field $\bar{\bar{\theta}} = 300$ K and corresponding pressure and density fields in hydrostatic balance\footnote{This corresponds to a state such that the gravitational and vertical pressure gradient forces are in balance, i.e. $g + c_p \theta \partial_z \pi = 0$.}. Additionally, we assume a zero velocity background state. A temperature perturbation of the form
\begin{align}
\Delta T =
\begin{cases}
&-\frac{15}{2}(1 + \cos(r\pi)) \hspace{4mm} \text{if } r<1,\\
&0  \hspace{31mm} \text{otherwise},
\end{cases}
\hspace{2cm} r = \sqrt{\frac{(x - x_c)^2}{x_r^2} + \frac{(z - z_c)^2}{z_r^2}}, \label{radius}
\end{align}
is added to the background potential temperature, while the density field is left unperturbed. The perturbation's horizontal and vertical centre and radius are given by $(x_c, x_r) = (16, 4)$ and $(z_c, z_r) = (3, 2)$ kilometres, respectively. The gravitational acceleration is defined as in the thermal shallow water test case, and the remaining physical parameters are given by $c_v=716.5$ m$^2$s$^{-2}$K$^{-1}$, $R=287$ m$^2$s$^{-2}$K$^{-1}$, and $p_0=100$ kPa. Finally, we note that a constant viscosity term is commonly added to the continuous equations in this test case to obtain a solution that converges as the resolution is refined \cite{melvin2010inherently}. In our case, we do not include this term due to our focus on energy conservation.\\ \\
The mesh is given by a vertically extruded interval mesh, with horizontal and vertical resolutions equal to $\Delta x = \Delta z = 100$ m. The simulation is run for 900 seconds, with a time step of $\Delta t = 0.5$ s, and 32 Picard iterations for each time step.\\ \\
To focus on the effects of upwind stabilisation, we consider the field development after 400 and 800 seconds. The resulting images and relative energy error  as well as $D\!G_\rho$ and $D\!G_\mathbf{u}$ developments are depicted in Figures \ref{DB_potential_temperatures}, \ref{DB_energies}, and \ref{DB_DG_rho_u}, in an arrangement equal to the one for the thermal shallow water equations above.\\
\\
\\
\\
We find that if the SUPG method is not applied in the vertical direction of potential temperature transport, spurious upward moving features occur in the potential temperature field development near the centre of the domain, where the bubble falls towards the bottom boundary. In the presence of the SUPG method these features are removed, again indicating a qualitatively favourable field development if the method is included. Further, as expected, energy is conserved up to solver tolerance for the two energy conserving brackets, with occasional jumps likely due to an insufficient number of Picard iterations. In contrast, the non-energy conserving bracket leads to a loss of energy of the order of $10^{-7}$, 4 orders of magnitude larger than the energy conserving ones. Finally, as for the perturbed steady state test case considered for the thermal rotating shallow water equations, we find that the energy-conserving setup without SUPG leads to increased values of $DG_\mathbf{u}$. Again, this likely due to thermal field noise entering the velocity field via the pressure gradient term. For $DG_\rho$, we find a slightly different field development for the three schemes once the bubble enters a more turbulent regime as it moves along the bottom boundary. Altogether, as before the energy conserving scheme including SUPG does not lead to an increase in oscillations in the velocity and density fields.
\subsubsection{Fully upwind stabilised, approximately energy conserving scheme}
\begin{figure}[ht]
\begin{center}
\includegraphics[width=1\textwidth]{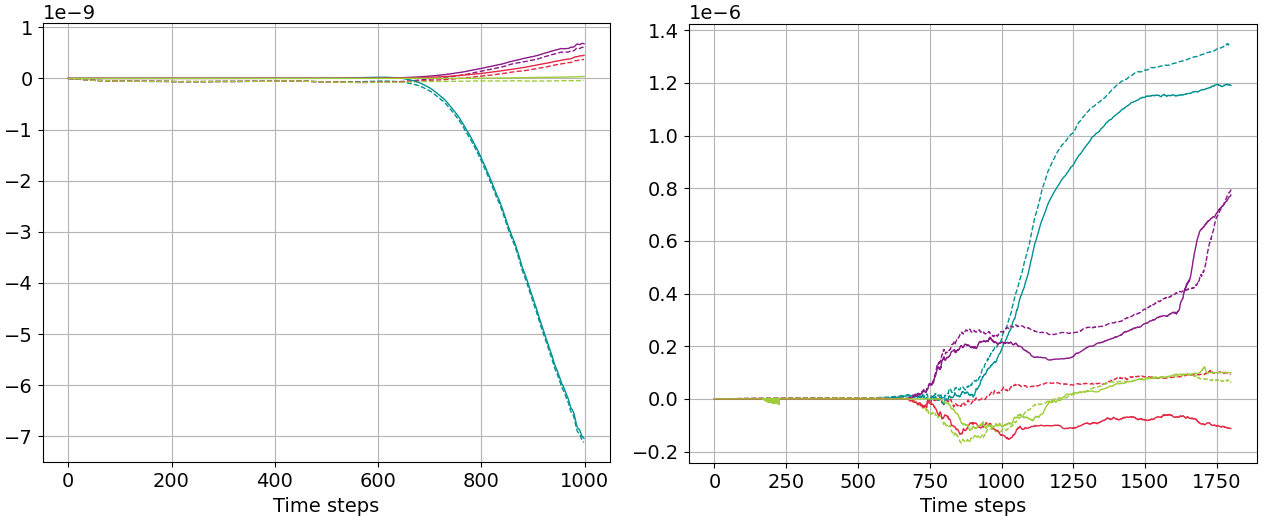}
\caption{Relative energy error development for rising (left) and falling (right) bubble test cases. Dashed lines indicate the approximated Poisson integrator, solid ones the full Poisson integrator. Colours correspond to cyan for 4 Picard iterations, purple for 5, red for 6, and light green for 8.} \label{nk_Poisson_vs_midpoint}
\end{center}
\end{figure}

\begin{figure}[ht]
\begin{center}
\includegraphics[width=1\textwidth]{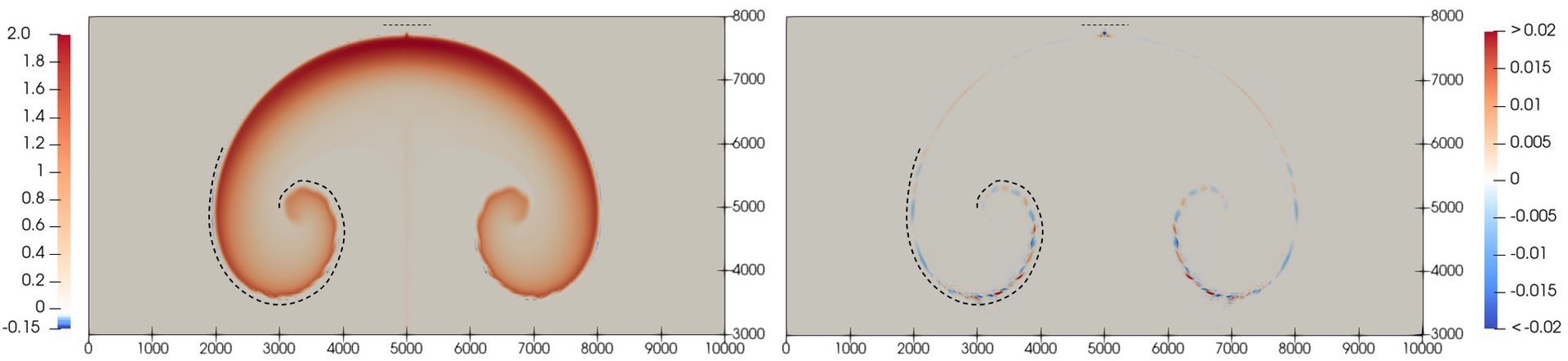}
\caption{Rising bubble test case at $1000$ s. Left: Potential temperature field for approximately energy conserving scheme. Right: potential temperature field difference between approximately energy conserving and non-energy conserving schemes. Dashed lines have been included to indicate relative positions between left and right plots. $\Delta x = \Delta z = 50$ m, $\Delta t = 0.25$ s.} \label{RB_approx}
\end{center}
\end{figure}

\begin{figure}[ht]
\begin{center}
\includegraphics[width=1\textwidth]{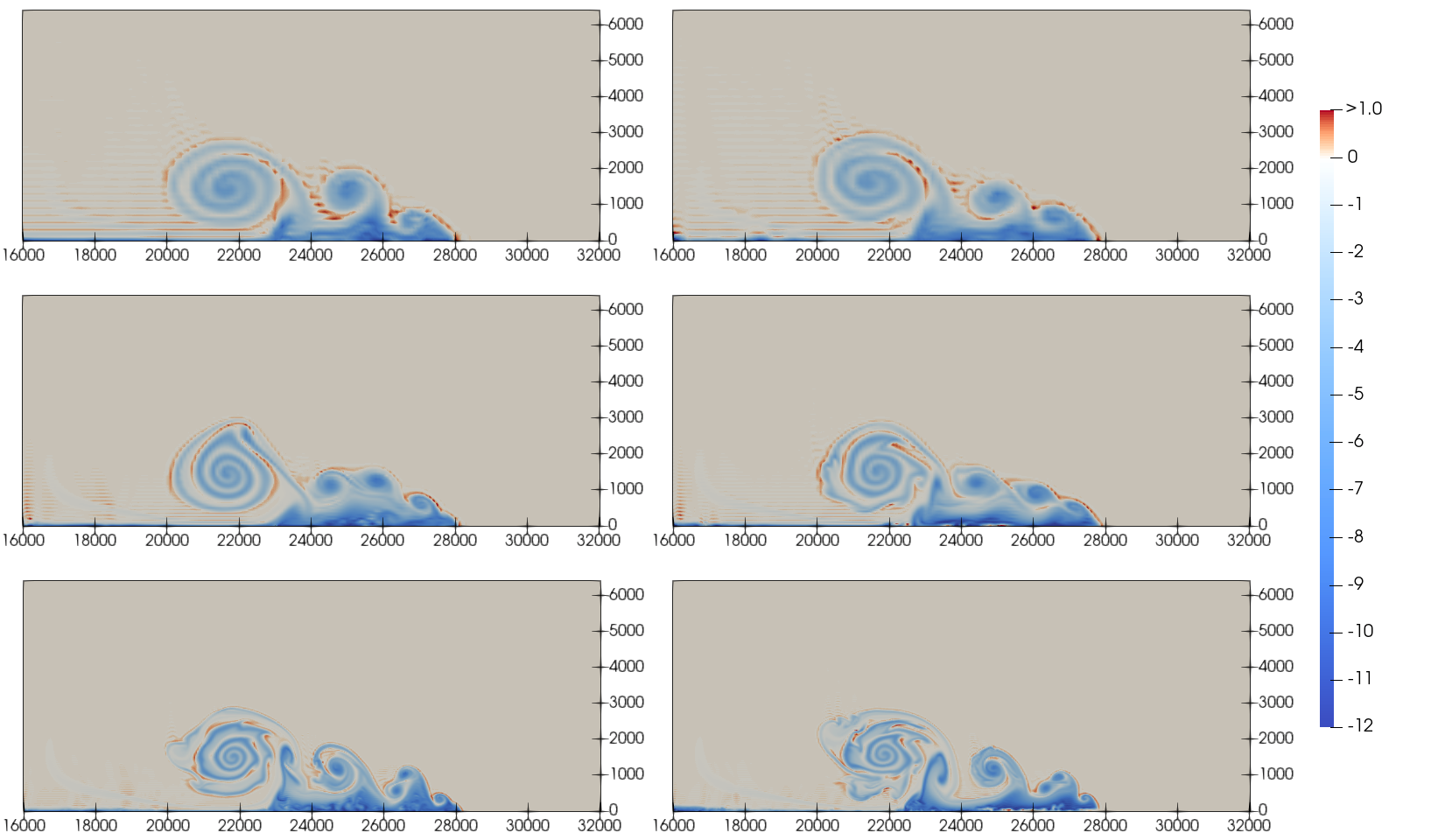}
\caption{Potential temperature fields for approximately energy conserving scheme (right column) and non-energy conserving scheme (left column), for falling bubble test case at 800 s. Top to bottom correspond to resolutions $200$ m, $100$ m, and $50$ m.} \label{DB_approx}
\end{center}
\end{figure}

\begin{figure}[ht]
\begin{center}
\includegraphics[width=1\textwidth]{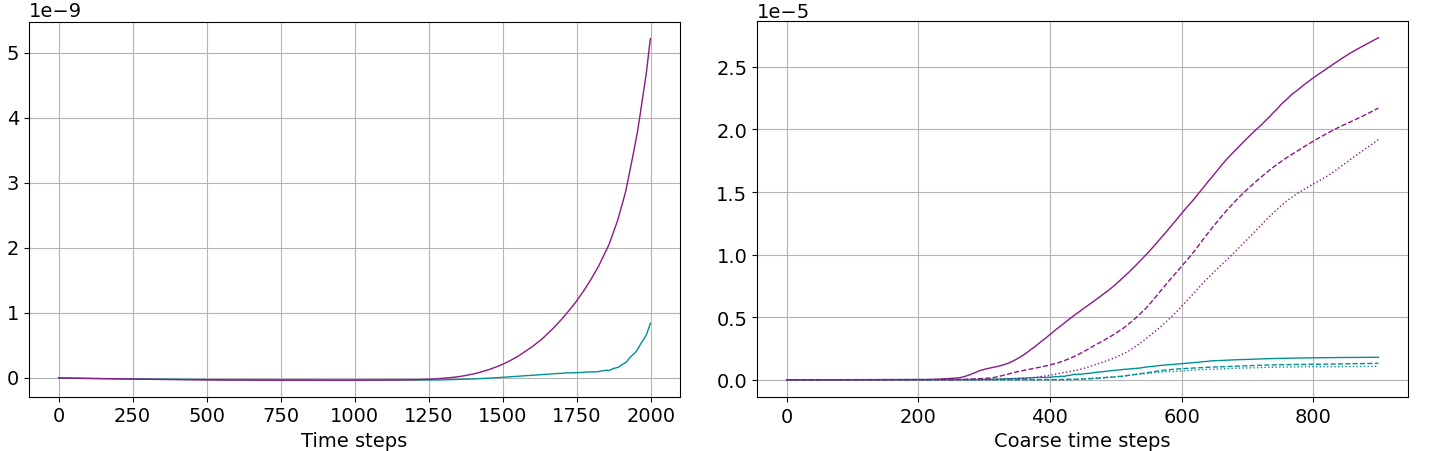}
\caption{Relative energy error development for approximately energy conserving scheme (cyan) and non-energy conserving scheme (purple), for 4 Picard iterations per time step. Left: rising bubble test case. Right: falling bubble test case, for resolutions $200$ m (solid), $100$ m (dashed), and $50$ m (dotted), and time step $\Delta t = 1$ s corresponding to $200$ m spatial resolution.} \label{full_upw_energies}
\end{center}
\end{figure}

While including the SUPG method for thermal field transport leads to an improvement of the field development, the latter still suffers from an absence of upwind stabilisation in the other fields, leading to an insufficient resolution of the density current flowing along the bottom boundary in the Euler test case. Further, 32 Picard iterations were required to achieve energy conservation of the order of $10^{-11}$. We therefore next consider the field and energy error developments for the approximately energy conserving scheme as presented in Section \ref{section_approx_ec}, which is energy conserving in space and energy conserving in time up to the difference given by the discretisation of the Hamiltonian variation in $\mathbf{u}$. The test cases considered here are given by the falling bubble one as described above, as well as a rising bubble test case based on \cite{carpenter1990application}. The latter is prone to secondary plumes (see e.g. \cite{bendall2019compatible}), and to avoid these in our discussion here, we consider the next higher order set of finite element spaces (k=3) for this test case.
The domain $\Omega$ is given by a horizontally periodic square of $10$ km side length. The background fields are given as in the falling bubble case, and the initial potential temperature perturbation is given by
\begin{align}
\Delta \theta =
\begin{cases}
&2 \cos^2(\frac{r\pi}{2}) \hspace{4mm} \text{if } r<1,\\
&0  \hspace{18mm} \text{otherwise},
\end{cases}
\end{align}
for $r$ as in \eqref{radius} above, with $(x_c, x_r) = (5, 2)$ and $(z_c, z_r) = (2, 2)$ kilometres. As before, the density field is left unperturbed. For the relative energy error study below, the mesh is as for the falling bubble test case with $\Delta x = \Delta z = 100$ m. For the rising bubble field development figures below, a higher resolution of $\Delta x = \Delta z = 50$ m is used. The simulations are run for 1000 seconds, with a time step of $\Delta t = 1$ s for the 100 m mesh and $\Delta t = 0.5$ s for the 50 m mesh.\\ \\
First, we compare the impact of both the reduced number of Picard iterations and the approximated Poisson integrator (i.e. \eqref{midpoint_dHdu}) on the relative energy error development. For this purpose, we run the falling and rising bubble test cases for $k \in \{4, 5, 6, 8\}$ Picard iterations, each time for both the fully energy conserving Poisson integrator and the approximated version as given in Section \ref{section_approx_ec}. The resulting energy developments are depicted in Figure \ref{nk_Poisson_vs_midpoint}. \\
\\
While more evident for the less turbulent rising bubble test case, we find that in both cases the choice of number of Picard iterations dominates the relative energy error development when compared to whether or not the Poisson integrator has been approximated. In particular, this indicates that the approximated scheme as given in Section \ref{section_approx_ec} does not lead to a substantial increase in energy error when a small number of Picard iterations is used, suggesting that the approximation (and therefore the removal of the additional computation for $\mathbb{U}$ appearing in the fully upwinded, fully energy conserving discretisation) is justified.\\ \\
To test the improvement of the qualitative field development when a fully upwinded bracket is used, we compare the scheme as given in Section \ref{section_approx_ec} with a fully upwinded reference scheme derived in a non-Hamiltonian setup (i.e. directly from the continuous equations), as defined in \eqref{u_eqn_nec_mid} - \eqref{theta_eqn_nec_mid}. We also consider the energy development for this comparison, using 4 Picard iterations in each time step. Note that considering the two discretisations, we find that the only additional computational cost of the approximately energy conserving scheme lies in computing $\overline{\dd{H}{\rho}}$ and $s(\tau, \bar{\mathbf{u}}, \bar{T})$ for each Picard iteration, i.e. a DG space and a horizontally DG space projection. To additionally study the effect of energy conservation depending on the choice of resolution, we consider the falling bubble test case at resolutions $\Delta x = \Delta z = 50 \text{ m}, 100 \text{ m}, 200 \text{ m}$, with time steps $\Delta t = 0.25 \text{ s}, 0.5 \text{ s}, 1\text{ s}$, respectively.\\ \\
For the rising bubble test case, the resulting potential temperature fields of the two schemes at $t = 1000$ seconds are visually difficult to distinguish, and we therefore consider in Figure \ref{RB_approx} the field development of the approximately energy conserving scheme, together with a plot depicting the difference between the latter scheme and the non-energy conserving one. In contrast, differences can clearly be seen for the falling bubble test case, and Figure \ref{DB_approx} contains the potential temperature fields for both schemes at $t = 800$ seconds. Finally, the corresponding relative energy error developments are given in Figure \ref{full_upw_energies}.
\\
\\
\\
For the rising bubble test case, the potential temperature field development of the two schemes differs only by a relatively small range compared to the initial perturbation magnitude $\Delta \theta \approx 2 K$. Considering the difference plot, we find that the largest discrepancies occur at regions where small scale features are forming. The first is given by the secondary plume at the bubble's top (with values differing by up to $\sim\! 0.12$ K). The second region is given by the bubble's front along the trailing vortices, where shear vortices are beginning to form (with values differing by up to $\sim\! 0.05$ K). Note that for the purpose of a better visibility of the field difference, a smaller range than these largest difference values has been chosen in Figure \ref{RB_approx}, and the corresponding secondary plume and shear vortex regions are shaded in dark red/dark blue. Finally, in terms of the relative energy error, we find that the approximately energy conserving scheme leads to a relative energy error that is reduced by a factor of $\sim\! 6.5$ when compared to the non-energy conserving scheme.\\ \\
Similarly, for the falling bubble test case we find an improved energy conservation by nearly an order of magnitude. Note that should a higher degree of energy conservation be required (albeit at a higher computational cost),  this can readily be achieved in the approximately energy conserving setup by using a moderately higher number of Picard iterations (with another order of magnitude gained for a total of 6 iterations per time step, as indicated in Figure \ref{nk_Poisson_vs_midpoint}). Additionally, we observe a qualitatively different field development depending on the resolution, which -- as indicated by the rising bubble test case -- is related to small scale features. The approximately energy conserving discretisation exhibits such features in the form of shear vortices within the main vortex at a resolution of $100$ m, while this only occurs at the next finer resolution for the non-Hamiltonian discretisation. This demonstrates that even in simulations with a comparable computational effort, the approximately energy conserving discretisation leads to a significant reduction in the relative energy error in comparison to a corresponding non-Hamiltonian discretisation, which in turn may lead to small scale features appearing at relatively coarser resolutions.
\begin{remark} \label{sub_grid_features}
The additional smaller vortices in the convergence study depicted in Figure \ref{DB_approx} also occur in a similar, viscosity-free falling bubble test case in \cite{marras2015stabilized} (see Figures 4 to 7), where a dynamic sub-grid scale model is used to stabilise a higher order spectral element discretisation. As discussed above, these features depend on the choice of resolution. In our case, the features' appearance at varyingly fine resolutions can be interpreted as a consequence of the higher degree of energy conservation, which in turn follows from the antisymmetric structure of the underlying almost Poisson bracket. As commented in Section \ref{section_SUPG_bracket}, the antisymmetry ensures correct transfers between kinetic, potential, and internal energy, in the sense that the growths and decreases of these three types of energy will not match exactly if the bracket is not antisymmetric. Indeed, given the $100$ m resolution setup, as we move from the bracket that underlies the approximately energy conserving setup \eqref{u_eqn_mid} - \eqref{theta_eqn_mid} term by term towards the non-antisymmetric formulation \eqref{u_eqn_nec_mid} - \eqref{theta_eqn_nec_mid}, the features become less pronounced and eventually disappear (figures not shown here). In ongoing work, we are exploring the energy transfer mechanisms leading to this behaviour.
\end{remark}
\section{Conclusion} \label{section_conclusion}
In this paper, we presented an almost Poisson bracket discretisation using the compatible finite element method, which includes an SUPG method for the thermal field transport equation. The bracket was used to derive fully energy conserving schemes for the Euler and thermal rotating shallow water equations. The SUPG formulation relies on an additional operator that corresponds to a SUPG-modified mass matrix. We demonstrated how to recover the thermal and velocity field evolution equations given the SUPG-modified bracket, showing that the thermal field equation corresponds to the usual SUPG formulation of a transport equation with the given underlying finite element spaces. Further, for the interest of a fully upwind stabilised, computationally less costly discretisation, we introduced an approximately energy conserving scheme including upwind stabilisation in the density and velocity fields. \\ \\
In numerical tests, we demonstrated energy conservation up to solver tolerance for the SUPG-modified bracket for an Euler and a thermal rotating shallow water scenario. Further, in both cases the incorporation of the SUPG method was shown to lead to a significant improvement of the qualitative thermal field development. Finally, we considered the approximately energy conserving discretisation in a simulation with a small number of Picard iterations, showing that we still gain a reduced energy error when compared to a reference discretisation derived in a non-Hamiltonian setup, at a comparable computational cost. For the more turbulent of the two Euler test cases considered here, the reduced energy error additionally leads to the appearance of small scale features starting from a coarser resolution than the reference discretisation.\\ \\
In future work, we aim to compare different upwind stabilisation methods within the Hamiltonian framework. In particular, this includes the choice of numerical method for velocity transport, which can e.g. be either formulated using upwind stabilisation for velocity, or a SUPG-type formulation for vorticity as presented for the shallow water equations in \cite{BAUER2018171}. While the numerical tests for the SUPG method considered here were restricted to two dimensions, this comparison will be done in a fully three-dimensional setting. We also aim to further investigate the role of balanced transfers between kinetic, potential, and internal energies, which -- next to total energy conservation -- are guaranteed by the Poisson bracket framework.\\ \\ \\
\textbf{Acknowledgements}\\
GAW acknowledges support from the EPSRC Mathematics of Planet Earth Centre for Doctoral Training at Imperial College London and the  University of Reading. CJC and WB acknowledge support from EPSRC grant EP/R029628/1 and NERC grant NE/M013634/1. WB also acknowledges the support of the `Laboratoires d'Excellence' CominLabs, Lebesgue and Mer through the SEACS project. Further, the authors are grateful to the anonymous reviewer who proposed a simplification to the original manuscript that led to the formulation using the modified mass matrix operator $s$.
\appendix
\renewcommand{\thesection}{\Alph{section}}
\section{Invertibility of SUPG-modified mass matrix} \label{appendix_A}
In this section, we show that $s(\tau, \mathbf{u}; \gamma)$ is well-defined. For this purpose, we first note that its definition \eqref{def_s_form} can be formulated in terms of a bilinear form $M_s$ defined on the thermal field space $\mathbb{V}_\theta$ such that
\begin{align}
M_s(\sigma_1, \sigma_2) =  \left\langle \sigma_1 + \tau S(\mathbf{u}; \sigma_1), \sigma_2 \right\rangle && \forall \sigma_1, \sigma_2 \in \mathbb{V}_\theta, \label{def_Ms}
\end{align}
where the SUPG contribution $S$ is given by \eqref{def_upwind_gamma}. In particular, $s$ is then given by
\begin{align}
M_s(s, \sigma) =  \langle \gamma,  \sigma \rangle && \forall \sigma \in \mathbb{V}_\theta, \label{def_Ms_form}
\end{align}
and following the Lax-Milgram theorem (see e.g. \cite{brenner2007mathematical}), $s$ is well-defined if $M_s$ is coercive and continuous. For simplicity, we assume the stabilisation parameter $\tau$ to be constant, as done in the numerical results section. Further, we need to make assumptions on the velocity field, as detailed in the lemma below.
\begin{lemma} \label{M_s_cts_coercive}
Let $\mathbb{V}_u(\Omega)$ be a div-conforming finite element space defined on either a regular rectangular mesh $\mathcal{T}^h$ for a rectangular domain $\Omega$, or a quasi-uniform triangular mesh on a spherical domain $\Omega$. Further, let $\mathbb{V}_\theta(\Omega)$, also defined on $\mathcal{T}^h$, be a continuous Galerkin space if $\Omega$ is a spherical domain, and a Charney-Phillips type finite element space otherwise. Consider a velocity field $\mathbf{u}(\mathbf{x}, t) \in L^2\big([0, T]; \mathbb{V}_u(\Omega)\big)$, for some $T > 0$. Then if there exists a constant $c_0$ such that
\begin{equation}
\|\mathbf{u}\|_{\infty, \Omega \times [0, T]} < c_0,
\end{equation}
the bilinear form $M_s$ defined by \eqref{def_Ms} satisfies
\begin{align}
M_s(\sigma_1, \sigma_2) \le \left(1 + C \frac{\tau}{h}\right) \|\sigma_1\|_2\|\sigma_2\|_2 && \forall \sigma_1, \sigma_2 \in \mathbb{V}_\theta, \label{M_s_cts}
\end{align}
where $C$ depends on the mesh geometry and $c_0$. Further, if there exists a constant $c_1$ such that
\begin{equation}
\|\nabla \cdot \mathbf{u}\|_{\infty, \Omega \times [0, T]} < c_1 \;\;\; (\mathbb{V}_\theta \text{ CG space}) \;\;\; \text{or} \;\;\; \|\partial_z u_z\|_{\infty, \Omega \times [0, T]} < c_1 \;\;\; (\mathbb{V}_\theta \text{ CP space}), \label{norm_u_der}
\end{equation}
where $z$ denotes the rectangle's second coordinate and $u_z$ the second component of $\mathbf{u}$, then $M_s$ additionally satisfies
\begin{align}
M_s(\sigma, \sigma) \ge \left(1 - \frac{c_1\tau}{2}\right) \|\sigma\|_2^2 && \forall \sigma \in \mathbb{V}_\theta. \label{M_coercive}
\end{align}
\textbf{Proof.} First, note that $M_s$ consists of a mass and an SUPG-type term, i.e.
\begin{align}
M_s(\sigma_1, \sigma_2) = \langle \sigma_1, \sigma_2 \rangle + \left\langle \tau S(\mathbf{u}; \sigma_1), \sigma_2 \right\rangle && \forall \sigma \in \mathbb{V}_\theta.
\end{align}
If $\mathbb{V}_\theta$ is a CG space, we bound the second term according to
\begingroup
\addtolength{\jot}{2mm}
\begin{align}
\langle \tau \mathbf{u} \cdot \nabla \sigma_1, \sigma_2 \rangle &\le \tau \|\mathbf{u}\|_{\infty, \Omega\times[0, T]}\|\nabla \sigma_1\|_2\|\sigma_2\|_2 \nonumber \\
&\le \tau \|\mathbf{u}\|_{\infty, \Omega\times[0, T]}(\|\sigma_1\|_2^2 + \|\nabla \sigma_1\|_2^2)^{\frac{1}{2}}\|\sigma_2\|_2 \nonumber \\
&\le \tau c_0\|\sigma_1\|_{W_2^1(\Omega)}\|\sigma_2\|_2 \nonumber \\
&\le C\frac{\tau}{h} \|\sigma_1\|_2\|\sigma_2\|_2 & \forall \sigma_1, \sigma_2 \in \mathbb{V}_\theta,
\end{align}
\endgroup
where we applied the Cauchy-Schwarz inequality in the first line, and further used an inverse inequality to bound the $W_2^1(\Omega)$ norm by the $L^2-norm$ (see e.g. \cite{brenner2007mathematical}), noting that the latter step introduces a mesh geometry dependent constant and the resolution dependent factor $1/h$. The desired estimate \eqref{M_s_cts} then follows by combining the above with the mass term. The derivation is analogous for the Charney-Phillips setup, noting that in the SUPG expression $S$, $\mathbf{k} \cdot \mathbf{u} = u_z$ and
\begin{equation}
\|u_z\|_{\infty, \Omega\times[0, T]} < \|\mathbf{u}\|_{\infty, \Omega\times[0, T]}.
\end{equation}
Further, in this case we also use an adjusted version of the $W^{1,2}(\Omega)$-norm, where $\|\nabla \sigma_1\|_2^2$ is replaced by $\|\partial_z \sigma_1\|_2^2$. Unlike the $W^{1,2}(\Omega)$-norm, this norm is valid here since $\sigma_1$ is continuous in the $z$-direction.\\ \\
To demonstrate \eqref{M_coercive}, we again consider the mass and SUPG terms of $M_s$ separately, and if $\mathbb{V}_\theta$ is a CG space, we find
\begin{align}
\langle \tau \mathbf{u} \cdot \nabla \sigma, \sigma \rangle = \left\langle \tau \mathbf{u}, \frac{1}{2}\nabla \sigma^2 \right\rangle&= - \frac{\tau}{2} \left( \langle \nabla \cdot \mathbf{u}, \sigma^2 \rangle - \int_\Gamma \left[\!\left[(\mathbf{n}\cdot \mathbf{u})\sigma^2\right]\!\right]dS \right) \nonumber \\
&\ge -\frac{\tau}{2} \|\nabla \cdot \mathbf{u}\|_{\infty, \Omega \times [0, T]} \|\sigma\|_2^2 \ge -\frac{c_1\tau}{2} \|\sigma\|_2^2 && \forall \sigma \in \mathbb{V}_\theta,
\end{align}
where we used integration by parts. Note that the integral over all facets $\Gamma$ vanishes since $\sigma$ and $\mathbf{n} \cdot \mathbf{u}$ (for facet normal vector $\mathbf{n}$) are continuous across facets. Together with the mass term, this leads to the desired inequality \eqref{M_coercive}. An analogous derivation holds for the Charney-Phillips setup, noting that in this case the facet integral is given by
\begin{equation}
\int_{\Gamma_h}  \left[\!\left[(\mathbf{n}\cdot \mathbf{k})(\mathbf{k} \cdot \mathbf{u})\sigma^2\right]\!\right]dS = \int_{\Gamma_h}  \left[\!\left[(n_z u_z)\sigma^2\right]\!\right]dS, \label{facet_term_CP}
\end{equation}
where $\Gamma_h$ denotes the set of all horizontal interior facets, and $n_z$ is given by the component of $\mathbf{n}$ in the $z$-direction. Note that $n_z u_z$ is continuous across horizontal facets since we assume a regular rectangular mesh for the Charney-Phillips setup. \hfill $\square$
\end{lemma}
The lemma immediately implies the proposition below, noting that coercivity requires a choice of stabilization parameter such that $\tau \le 2/c_1$. Further, note that the continuity bound can be made independent of $h$ if $\tau$ is of the order of $h$. The coercivity criterion depends on the size of the $\mathbf{u}$-related terms \eqref{norm_u_der}, and in practice, we found $\tau = \Delta t/2$ to be sufficiently small. This choice of $\tau$ is also compatible with mesh independence, noting that $\tau$ can be written as
\begin{equation}
\tau = h \; \frac{c}{2|\mathbf{u}|},
\end{equation}
for Courant number $c$.\\ \\
Further, for \eqref{facet_term_CP} to vanish, we assumed that all horizontal interior facets of the underlying mesh are orthogonal to the $z$-coordinate's unit vector $\mathbf{k}$. This assumption will generally not hold true if we consider near rectangular domains with meshes not aligned with the $x$-coordinate. In particular, the latter is the case for atmospheric domains where topography is included and a terrain-following mesh is used, in which case further care needs to be taken to bound the facet integral accordingly. Finally, the simplifying assumption of constant $\tau$ can also be relaxed (see e.g. \cite{bochev2004stability} for coercivity in the case of divergence-free $\mathbf{u}$). However, since this paper is focused on incorporating the SUPG method in an energy conserving framework, we deliberately kept the discussion on the choice of $\tau$ parameter simple.
\begin{prop}
Assume the velocity field $\mathbf{u}$ satisfies the conditions given in Lemma \ref{M_s_cts_coercive}, and further that the SUPG parameter $\tau$ is such that $\tau \le 2/c_1$, for $c_1$ as given in \eqref{norm_u_der}. Then the SUPG operator $s$ given in Definition \ref{definition_s} is well-defined and invertible.\\ \\
\text{Proof.} Follows directly from the Lax Milgram theorem and Lemma \ref{M_s_cts_coercive}. \hfill $\square$
\end{prop}
\newpage
\bibliography{Poisson_SUPG}

\begin{thebibliography}{10}

\bibitem{arnold2002unified}
D.~N. Arnold, F.~Brezzi, B.~Cockburn, and L.~D. Marini.
\newblock Unified analysis of discontinuous {G}alerkin methods for elliptic
  problems.
\newblock {\em SIAM journal on numerical analysis}, 39(5):1749--1779, 2002.

\bibitem{balay2019petsc}
S.~Balay, S.~Abhyankar, M.~Adams, J.~Brown, P.~Brune, K.~Buschelman, L.~Dalcin,
  A.~Dener, V.~Eijkhout, W.~D. Gropp, et~al.
\newblock {P}{E}{T}{S}c users manual.
\newblock 2019.

\bibitem{balay1997efficient}
S.~Balay, W.~D. Gropp, L.~C. McInnes, and B.~F. Smith.
\newblock Efficient management of parallelism in object-oriented numerical
  software libraries.
\newblock In {\em Modern software tools for scientific computing}, pages
  163--202. Springer, 1997.

\bibitem{BAUER2018171}
W.~Bauer and C.~J. Cotter.
\newblock Energy--enstrophy conserving compatible finite element schemes for
  the rotating shallow water equations with slip boundary conditions.
\newblock {\em Journal of Computational Physics}, 373:171 -- 187, 2018.

\bibitem{bauer2019variational}
W.~Bauer and F.~Gay-Balmaz.
\newblock Towards a geometric variational discretization of compressible
  fluids: the rotating shallow water equations.
\newblock {\em Journal of Computational Dynamics}, 6(2158-2491):1, 2019 2019.

\bibitem{bendall2019compatible}
T.~M. Bendall, T.~H. Gibson, J.~Shipton, C.~J. Cotter, and B.~Shipway.
\newblock A compatible finite element discretisation for the moist compressible
  {E}uler equations.
\newblock {\em arXiv preprint arXiv:1910.01857}, 2019.

\bibitem{gmd-9-3803-2016}
G.-T. Bercea, A.~T.~T. McRae, D.~A. Ham, L.~Mitchell, F.~Rathgeber, L.~Nardi,
  F.~Luporini, and P.~H.~J. Kelly.
\newblock A structure-exploiting numbering algorithm for finite elements on
  extruded meshes, and its performance evaluation in {F}iredrake.
\newblock {\em Geoscientific Model Development}, 9(10):3803--3815, 2016.

\bibitem{bochev2004stability}
P.~B. Bochev, M.~D. Gunzburger, and J.~N. Shadid.
\newblock Stability of the {SUPG} finite element method for transient
  advection--diffusion problems.
\newblock {\em Computer methods in applied mechanics and engineering},
  193(23-26):2301--2323, 2004.

\bibitem{brenner2007mathematical}
S.~Brenner and R.~Scott.
\newblock {\em The mathematical theory of finite element methods}, volume~15.
\newblock Springer Science \& Business Media, 2007.

\bibitem{brooks1982streamline}
A.~N. Brooks and T.~J.~R. Hughes.
\newblock Streamline upwind/{P}etrov-{G}alerkin formulations for convection
  dominated flows with particular emphasis on the incompressible
  {N}avier-{S}tokes equations.
\newblock {\em Computer Methods in Applied Mechanics and Engineering},
  32(1-3):199--259, 1982.

\bibitem{carpenter1990application}
R.~L. Carpenter~Jr, K.~K. Droegemeier, P.~R. Woodward, and C.~E. Hane.
\newblock Application of the piecewise parabolic method ({P}{P}{M}) to
  meteorological modeling.
\newblock {\em Monthly Weather Review}, 118(3):586--612, 1990.

\bibitem{cohen2011linear}
D.~Cohen and E.~Hairer.
\newblock Linear energy-preserving integrators for {P}oisson systems.
\newblock {\em BIT Numerical Mathematics}, 51(1):91--101, 2011.

\bibitem{COTTER20127076}
C.~J. Cotter and J.~Shipton.
\newblock Mixed finite elements for numerical weather prediction.
\newblock {\em Journal of Computational Physics}, 231(21):7076 -- 7091, 2012.

\bibitem{ELDRED20191}
C.~Eldred, T.~Dubos, and E.~Kritsikis.
\newblock A quasi-{H}amiltonian discretization of the thermal shallow water
  equations.
\newblock {\em Journal of Computational Physics}, 379:1 -- 31, 2019.

\bibitem{ford2013gung}
R.~Ford, M.~J. Glover, D.~A. Ham, C.~M. Maynard, S.~M. Pickles, G.~Riley, and
  N.~Wood.
\newblock {G}ung {H}o: A code design for weather and climate prediction on
  exascale machines.
\newblock In {\em Proceedings of the Exascale Applications and Software
  Conference}, 2013.

\bibitem{gassmann2013global}
A.~Gassmann.
\newblock A global hexagonal {C}-grid non-hydrostatic dynamical core
  ({I}{C}{O}{N}-{I}{A}{P}) designed for energetic consistency.
\newblock {\em Quarterly Journal of the Royal Meteorological Society},
  139(670):152--175, 2013.

\bibitem{gibson2020slate}
T.~H. Gibson, L.~Mitchell, D.~A. Ham, and C.~J. Cotter.
\newblock Slate: extending {F}iredrake's domain-specific abstraction to
  hybridized solvers for geoscience and beyond.
\newblock {\em Geoscientific model development}, 13(2):735--761, 2020.

\bibitem{kirby2018solver}
R.~C. Kirby and L.~Mitchell.
\newblock Solver composition across the {P}{D}{E}/linear algebra barrier.
\newblock {\em SIAM Journal on Scientific Computing}, 40(1):C76--C98, 2018.

\bibitem{kuzmin2010guide}
D.~Kuzmin.
\newblock A guide to numerical methods for transport equations.
\newblock {\em University Erlangen-Nuremberg}, 2010.

\bibitem{lee2021petrov}
D.~Lee.
\newblock {P}etrov--{G}alerkin flux upwinding for mixed mimetic spectral
  elements, and its application to geophysical flow problems.
\newblock {\em Computers \& Mathematics with Applications}, 89:68--77, 2021.

\bibitem{lee2020mixed}
D.~Lee and A.~Palha.
\newblock A mixed mimetic spectral element model of the 3{D} compressible
  {E}uler equations on the cubed sphere.
\newblock {\em Journal of Computational Physics}, 401:108993, 2020.

\bibitem{lucarini2011energetics}
V.~Lucarini and F.~Ragone.
\newblock Energetics of climate models: net energy balance and meridional
  enthalpy transport.
\newblock {\em Reviews of Geophysics}, 49(1), 2011.

\bibitem{marras2015stabilized}
S.~Marras, M.~Nazarov, and F.~X. Giraldo.
\newblock Stabilized high-order {G}alerkin methods based on a parameter-free
  dynamic {SGS} model for {LES}.
\newblock {\em Journal of Computational Physics}, 301:77--101, 2015.

\bibitem{mcrae2016automated}
A.~T.~T. McRae, G.-T. Bercea, L.~Mitchell, D.~A. Ham, and C.~J. Cotter.
\newblock Automated generation and symbolic manipulation of tensor product
  finite elements.
\newblock {\em SIAM Journal on Scientific Computing}, 38(5):S25--S47, 2016.

\bibitem{melvin2018choice}
T.~Melvin, T.~Benacchio, J.~Thuburn, and C.~J. Cotter.
\newblock Choice of function spaces for thermodynamic variables in mixed
  finite-element methods.
\newblock {\em Quarterly Journal of the Royal Meteorological Society},
  144(712):900--916, 2018.

\bibitem{melvin2010inherently}
T.~Melvin, M.~Dubal, N.~Wood, A.~Staniforth, and M.~Zerroukat.
\newblock An inherently mass-conserving iterative semi-implicit
  semi-{L}agrangian discretization of the non-hydrostatic vertical-slice
  equations.
\newblock {\em Quarterly Journal of the Royal Meteorological Society},
  136(648):799--814, 2010.

\bibitem{morrison1998hamiltonian}
P.~J. Morrison.
\newblock {H}amiltonian description of the ideal fluid.
\newblock {\em Reviews of modern physics}, 70(2):467, 1998.

\bibitem{morrison1980noncanonical}
P.~J. Morrison and J.~M. Greene.
\newblock Noncanonical {H}amiltonian density formulation of hydrodynamics and
  ideal magnetohydrodynamics.
\newblock {\em Physical Review Letters}, 45(10):790, 1980.

\bibitem{natale2017scale}
A.~Natale and C.~J. Cotter.
\newblock Scale-selective dissipation in energy-conserving finite element
  schemes for two-dimensional turbulence.
\newblock {\em Quarterly Journal of the Royal Meteorological Society},
  143(705):1734--1745, 2017.

\bibitem{natale2016variational}
A.~Natale and C.~J. Cotter.
\newblock A variational {H}(div) finite element discretisation for perfect
  incompressible fluids.
\newblock {\em IMA Journal of Numerical Analysis}, 38(3):1388--1419, 2017.

\bibitem{natale2016compatible}
A.~Natale, J.~Shipton, and C.~J. Cotter.
\newblock Compatible finite element spaces for geophysical fluid dynamics.
\newblock {\em Dynamics and Statistics of the Climate System}, 1(1), 2016.

\bibitem{rathgeber2016firedrake}
F.~Rathgeber, D.~A. Ham, L.~Mitchell, M.~Lange, F.~Luporini, A.~T.~T. McRae,
  G.~T. Bercea, G.~R. Markall, and P.~H.~J. Kelly.
\newblock Firedrake: automating the finite element method by composing
  abstractions.
\newblock {\em ACM Transactions on Mathematical Software (TOMS)}, 43(3):24,
  2016.

\bibitem{ripa1993conservation}
P~Ripa.
\newblock Conservation laws for primitive equations models with inhomogeneous
  layers.
\newblock {\em Geophysical \& Astrophysical Fluid Dynamics}, 70(1-4):85--111,
  1993.

\bibitem{shepherd1990symmetries}
T.~G. Shepherd.
\newblock Symmetries, conservation laws, and {H}amiltonian structure in
  geophysical fluid dynamics.
\newblock {\em Adv. Geophys}, 32(287-338):2, 1990.

\bibitem{shipton2018higher}
J.~Shipton, T.~H. Gibson, and C.~J. Cotter.
\newblock Higher-order compatible finite element schemes for the nonlinear
  rotating shallow water equations on the sphere.
\newblock {\em Journal of Computational Physics}, 375:1121--1137, 2018.

\bibitem{straka1993numerical}
J.~M. Straka, R.~B. Wilhelmson, L.~J. Wicker, J.~R. Anderson, and K.~K.
  Droegemeier.
\newblock Numerical solutions of a non-linear density current: a benchmark
  solution and comparisons.
\newblock {\em International Journal for Numerical Methods in Fluids},
  17(1):1--22, 1993.

\bibitem{taylor2020energy}
M.~A. Taylor, O.~Guba, A.~Steyer, P.~A. Ullrich, D.~M. Hall, and C.~Eldred.
\newblock An energy consistent discretization of the nonhydrostatic equations
  in primitive variables.
\newblock {\em Journal of Advances in Modeling Earth Systems}, 12(1), 2020.

\bibitem{williamson1992standard}
D.~L. Williamson, J.~B. Drake, J.~J. Hack, R.~Jakob, and P.~N. Swarztrauber.
\newblock A standard test set for numerical approximations to the shallow water
  equations in spherical geometry.
\newblock {\em Journal of Computational Physics}, 102(1):211--224, 1992.

\bibitem{wimmer2020energy}
G.~A. Wimmer, C.~J. Cotter, and W.~Bauer.
\newblock Energy conserving upwinded compatible finite element schemes for the
  rotating shallow water equations.
\newblock {\em Journal of Computational Physics}, 401:109016, 2020.

\bibitem{zeitlin2018geophysical}
V.~Zeitlin.
\newblock {\em Geophysical fluid dynamics: understanding (almost) everything
  with rotating shallow water models}.
\newblock Oxford University Press, 2018.

\end{thebibliography}
\end{document}